\documentclass[preprint,3p,sort&compress,final,times]{elsarticle}
\usepackage{amsmath}
\usepackage{amssymb}
\usepackage{mathtools}
\usepackage{placeins}
\usepackage[usenames,dvipsnames]{color}
\usepackage{cases}
\usepackage{array}

\usepackage{lineno}
\usepackage[usenames,dvipsnames]{color}

\newcommand{\qMass}{\langle\gamma_h,\gamma_h\rangle}

\newcommand{\uMassThetaPar}{\langle\boldsymbol{\beta}_h^{\parallel},{\theta}_h^{n+1/2}\boldsymbol{\beta}_h^{\parallel}\rangle}
\newcommand{\uMassThetaPerp}{\langle\beta_h^{\perp},{\theta}_h^{n+1/2}\beta_h^{\perp}\rangle}
\newcommand{\uMassPar}{\langle\boldsymbol{\beta}_h^{\parallel},\boldsymbol{\beta}_h^{\parallel}\rangle}
\newcommand{\uMassPerp}{\langle\beta_h^{\perp},\beta_h^{\perp}\rangle}

\newcommand{\EdcPar}{\boldsymbol{\mathsf{E}}^{3,2}_{\parallel}}
\newcommand{\EdcPerp}{\boldsymbol{\mathsf{E}}^{3,2}_{\perp}}

\newcommand{\EdcTPar}{\Big(\boldsymbol{\mathsf{E}}^{3,2}_{\parallel}\Big)^{\top}}
\newcommand{\EdcTPerp}{\Big(\boldsymbol{\mathsf{E}}^{3,2}_{\perp}\Big)^{\top}}
\newcommand{\zero}{\boldsymbol{\mathsf{0}}}

\newcommand{\valgebra}[1]{\hat{#1}}

\newtheorem{remark}{Remark}

\newcommand{\blue}[1]{{#1}}
\newcommand{\green}[1]{{#1}}

\journal{}
\begin{document}
\begin{frontmatter}

\title{Exact spatial and temporal balance of energy exchanges within a
horizontally explicit/vertically implicit non-hydrostatic atmosphere}
\author[BOM]{David Lee\corref{cor}}
\ead{davelee2804@gmail.com}
\author[EIN]{Artur Palha}

\address[BOM]{Bureau of Meteorology, Melbourne, Australia}
\address[EIN]{Faculty of Aerospace Engineering, Delft University of Technology, Kluyverweg 2, 2629 HS Delft, The Netherlands}
\cortext[cor]{Corresponding author. Tel. +61 452 262 804.}

\begin{abstract}
A new horizontally explicit/vertically implicit (HEVI) time splitting scheme for atmospheric
modelling is introduced, for which the horizontal divergence terms are applied within the
implicit vertical substep. 
%When using a spatial and vertical temporal discretisation that
\blue{The new HEVI scheme is implemented in conjunction with a mixed mimetic spectral element 
spatial discretisation and semi-implicit vertical time stepping scheme that both }
preserve the skew-symmetric structure of the non-canonical Hamiltonian form of the equations
of motion. \blue{Within this context the new HEVI scheme } allows for the exact balance of 
all energetic exchanges in space and
time. However since the choice of horizontal fluxes for which this balance is satisfied is not
consistent with the horizontal velocity at the end of the time level the scheme still admits
a temporal energy conservation error. Linearised eigenvalue analysis shows that
similar to a fully implicit method, the new HEVI scheme is neutrally stable for all buoyancy
modes, and unlike a second order trapezoidal HEVI scheme is stable for all acoustic
modes below a certain horizontal CFL number. The scheme is validated against standard test
cases for both planetary and non-hydrostatic regimes. For the planetary scale baroclinic
instability test case, the new formulation exhibits a secondary oscillation in the potential
to kinetic energy power exchanges, with a temporal frequency approximately four times that
exhibited by a horizontally third order, vertically second order trapezoidal scheme. For the
non-hydrostatic test case, the vertical upwinding of the potential temperature diagnostic
equation is shown to reduce spurious oscillations without altering the energetics of the
solution, since this upwinding is performed in an energetically consistent manner.
\blue{For this test case, which is configured on an affine geometry, the exact balance of 
energy exchanges allows the model to run stably without any form of dissipation.}
\end{abstract}

\end{frontmatter}

Horizontally explicit/vertically implicit (HEVI) schemes are a popular approach to time
stepping in atmospheric models. This is on account of the importance of \blue{removing } the CFL
limit imposed by the explicit representation of vertical dynamics, and the desire to avoid
the computational expense of a three dimensional implicit solve at each nonlinear iteration.
Numerous HEVI schemes have been implemented \cite{Ullrich12,Giraldo13,Bao15,Gardner18,Steyer20}
in previous models. These schemes are typically motived by a desire to improve upon the
stability and dispersion properties of the second order trapezoidal HEVI scheme through the
use of additional sub steps within each time step \cite{Weller13,Lock14}. The application of
additional substeps qualitatively improves the stepping scheme by ensuring the linear stability
of the horizontal buoyancy and acoustic models below some CFL limit.

In the present article a new HEVI scheme is introduced, motivated by a desire to preserve
the exact balance of energetic exchanges. The scheme is implemented within the context
of a mixed mimetic spectral element spatial discretisation \cite{LP20} and energy conserving
implicit vertical time stepping method \cite{Lee20} that together allow for the exact balance
of all energy exchanges in space and time. This is achieved by evaluating the horizontal
mass and temperature fluxes within the vertically implicit sub-step, after first computing
a provisional horizontal velocity from which to derive a second order temporal representation
of the horizontal fluxes. Unlike other HEVI schemes, the new method is composed of only two
horizontal and one vertical substep. In contrast to the standard second order trapezoidal splitting scheme,
linear analysis shows the new method to be unconditionally stable for all acoustic modes
below some horizontal CFL number, and all buoyancy modes independent of CFL number. While the
energetic exchanges are all exactly balanced, since the provisional horizontal velocity does
not exactly match the final horizontal velocity at the end of the time step, the scheme admits
a horizontal energy conservation error, as is the case for explicit schemes used to
solve compressible flow problems. 

In section \ref{sec::ssvar} the continuous 3D compressible
Euler equations for atmospheric flows are introduced in variational form. The discrete analogue
of these equations \blue{and a discussion of the energy conserving time stepping are then presented } 
in section \ref{sec::df}, and the new HEVI splitting scheme is
described in section \ref{sec::hevi}. Section \ref{sec::upwind} describes the energetically consistent
stabilisation of spurious oscillations in potential temperature via the upwinding of the corresponding test space.
%Results for the linearised stability analysis are presented in section 4, and t
The reproduction of standard test cases is described in section \ref{sec::test}. Finally the
conclusions are discussed in section \ref{sec::conc}.

\section{Skew-symmetric variational structure of the compressible Euler equations}\label{sec::ssvar}

The compressible Euler equations for atmospheric flows may be expressed in skew-symmetric form
for the evolution of the velocity, $\boldsymbol{u}$, density, $\rho$, and density weighted
potential temperature, $\Theta$, as:
\begin{subequations}\label{eq::ce}
\begin{align}
\frac{\partial\boldsymbol{u}}{\partial t} &= -\boldsymbol{q}\times\boldsymbol{U} - \nabla\Phi - \theta\nabla\Pi\\
\frac{\partial\rho}{\partial t} &= -\nabla\cdot\boldsymbol{U}\label{eq::ce_rho}\\
\frac{\partial\Theta}{\partial t} &= -\nabla\cdot(\theta\boldsymbol{U})\label{eq::ce_theta}
\end{align}
\end{subequations}
where $\boldsymbol{q} := (\nabla\times\boldsymbol{u} + \boldsymbol{f})/\rho$ is the potential vorticity (and 
\green{$\boldsymbol{f} = f\boldsymbol{e}_z$ }
is the Coriolis term), $\boldsymbol{U} := \rho\boldsymbol{u}$ is the \blue{mass } flux,  $\Phi :=  \frac{1}{2}\boldsymbol{u}\cdot\boldsymbol{u} +  \rho g z$ is the Bernoulli function (with $g$ the gravitational acceleration and $z$ the height), $\theta := \Theta/\rho$ is the potential temperature, and $\Pi := c_{p}\left(R\Theta/p_{0}\right)^{\frac{R}{c_{v}}}$ is the Exner pressure (where $c_v$ and $c_p$ are the specific
heats at constant volume and temperature respectively, \green{$R$ is the ideal gas constant and $p_0$ is the reference pressure}). 
The corresponding Hamiltonian (total energy) is given as:
\begin{equation}
{H} = \int_{\Omega}\overbrace{\frac{1}{2}\rho\boldsymbol{u}\cdot\boldsymbol{u}}^{{K}} + \overbrace{\vphantom{\frac{1}{2}}\rho gz}^{{P}} + \overbrace{\frac{c_v}{c_p}\Theta\Pi}^{{I}}\mathrm{d}\Omega, \label{eq::total_energy}
\end{equation}
where ${K}$ is the kinetic energy, ${P}$ is the potential energy, and ${I}$ is the internal energy. 
\green{In this work the domain, $\Omega$ is configured as either the volume on the surface of a 
sphere in physical units, or a high resolution Cartesian geometry which is periodic in the two 
horizontal dimensions. } The variational derivatives of the
Hamiltonian with respect to the dependent variables ($\boldsymbol{u}$, $\rho$, $\Theta$) are
\begin{subequations} \label{eq::variational_derivatives_hamiltonian}
\begin{align}
\frac{\delta{H}}{\delta\boldsymbol{u}} &= \rho\boldsymbol{u} =: \boldsymbol{U}, \\
\frac{\delta{H}}{\delta\rho} &= \frac{1}{2}\boldsymbol{u}\cdot\boldsymbol{u} + gz =: \Phi, \\
\frac{\delta{H}}{\delta\Theta} &= c_p\Bigg(\frac{R\Theta}{p_0}\Bigg)^{\frac{R}{c_v}} =: \Pi ,
\end{align}
\end{subequations}
\blue{which } give the mass flux, $\boldsymbol{U}$, the Bernoulli function, $\Phi$ and the Exner pressure, $\Pi$ respectively.

The weak form of \eqref{eq::ce} is given for $\boldsymbol{u}, \boldsymbol{U}\in H(\mathrm{div},\Omega)$
and $\rho,\Theta,\Phi,\Pi\in L^2(\Omega)$ as:
\begin{subequations}\label{eq::ce_var}
\begin{align}
\Bigg\langle\boldsymbol{\beta},\frac{\partial\boldsymbol{u}}{\partial t}\Bigg\rangle &=
-\langle\boldsymbol{\beta},\boldsymbol{q}\times\boldsymbol{U}\rangle
+\langle\nabla\cdot\boldsymbol{\beta},\Phi\rangle
+\langle\nabla\cdot\mathbb{P}^{\mathrm{div}}\left[\theta\boldsymbol{\beta}\right], \Pi\rangle,\qquad\forall\boldsymbol{\beta}\in H(\mathrm{div},\Omega)\label{eq::ce_var_mom}\\
\Bigg\langle\gamma,\frac{\partial\rho}{\partial t}\Bigg\rangle &=
-\langle\gamma,\nabla\cdot\boldsymbol{U}\rangle,\qquad\forall\gamma\in L^2(\Omega)\\
\Bigg\langle\sigma,\frac{\partial\Theta}{\partial t}\Bigg\rangle &=
-\langle\sigma,\nabla\cdot\mathbb{P}^{\mathrm{div}}\left[\theta\boldsymbol{U}\right]\rangle,\qquad\forall\sigma\in L^2(\Omega)
\end{align}
\end{subequations}
where we have invoked inner product notation as $\langle a,b\rangle = \int_{\Omega} ab\mathrm{d}\Omega$, the integration-by-parts relation $\langle\boldsymbol{\beta},\nabla\alpha\rangle = -\langle\nabla\cdot\boldsymbol{\beta},\alpha\rangle$
(assuming periodic or homogeneous boundary conditions) in \eqref{eq::ce_var}, and introduced the $L^{2}$ projection to $H(\mathrm{div}, \Omega)$
	\begin{equation}
		\mathbb{P}^{\mathrm{div}}: \left[L^{2}(\Omega)\right]^{3} \rightarrow H(\mathrm{div}, \Omega)\,.
	\end{equation}
Specifically \green{for $\boldsymbol{v} \in \left[L^{2}(\Omega)\right]^{3}$, the } projection 
$\mathbb{P}^{\mathrm{div}}\left[\boldsymbol{v}\right]  \in H(\mathrm{div}, \Omega)$ is given by
\begin{equation}
	\langle\boldsymbol{\beta},\mathbb{P}^{\mathrm{div}}\left[\boldsymbol{v}\right]\rangle = \langle\boldsymbol{\beta},\boldsymbol{v}\rangle\,,\qquad \forall \boldsymbol{\beta}\in H(\mathrm{div}, \Omega)\,.
\end{equation}
Energy conservation is assured for the choice of $\boldsymbol{\beta} = \boldsymbol{U}$, $\gamma = \Phi$, $\sigma = \Pi$
by summing all equations in \eqref{eq::ce_var}, such that
\begin{equation}
\Bigg\langle\frac{\delta{H}}{\delta\boldsymbol{u}},\frac{\partial\boldsymbol{u}}{\partial t}\Bigg\rangle +
\Bigg\langle\frac{\delta{H}}{\delta\rho},\frac{\partial\rho}{\partial t}\Bigg\rangle +
\Bigg\langle\frac{\delta{H}}{\delta\Theta},\frac{\partial\Theta}{\partial t}\Bigg\rangle =
\frac{\mathrm{d}H}{\mathrm{d} t} = 0. \label{eq::hamiltonian_time_variation_0}
\end{equation}
where the key ingredient for energy conservation is the skew-symmetry of \eqref{eq::ce_var}.

\section{Discrete formulation}\label{sec::df}

In this section the mixed mimetic spectral element spatial discretisation used in the model is
introduced, \blue{as well as an analysis of the geometric properties of the energy conserving
time integration scheme}. 

\subsection{Spatial discretisation}

For the purposes of this article the salient feature of the \blue{mixed mimetic spectral element } method 
is that it allows for the preservation of the skew-symmetric structure of the compressible Euler equations,
and thus the conservation of energy and energetic exchanges in the discrete form.
For a more detailed discussion the reader is referred to previous work on the
use of this method for geophysical flow modelling \cite{LPG18,LP18,LP20}, as well as more foundational
works on the subject \cite{Gerritsma11,Kreeft13,Palha14,Hiemstra14}. In order to begin this discussion
we introduce the finite dimensional subspaces $\mathcal{P}_h\subset H^1(\Omega)$,
$\mathcal{W}_h\subset H(\mathrm{curl},\Omega)$, $\mathcal{U}_h\subset H(\mathrm{div},\Omega)$
and $\mathcal{Q}\subset L^2(\Omega)$. These subspaces are spanned by basis functions
\begin{equation}
	\mathrm{span} \{\psi_1, \dots, \psi_{\mathrm{d}_{\mathcal{P}}}\} = \mathcal{P}_h\,, \quad \mathrm{span} \{\boldsymbol{\alpha}_1, \dots, \boldsymbol{\alpha}_{\mathrm{d}_{\mathcal{W}}}\} = \mathcal{W}_h\,, \quad \mathrm{span} \{\boldsymbol{\beta}_1, \dots, \boldsymbol{\beta}_{\mathrm{d}_{\mathcal{U}}}\} = \mathcal{U}_h\,, \quad\text{and}\quad \mathrm{span} \{\gamma_1, \dots, \gamma_{\mathrm{d}_{\mathcal{Q}}}\} = \mathcal{Q}_h\,, \quad
\end{equation}
with $\mathrm{d}_{\mathcal{P}}$, $\mathrm{d}_{\mathcal{W}}$, $\mathrm{d}_{\mathcal{U}}$, and $\mathrm{d}_{\mathcal{Q}}$, the number of degrees of freedom of the spaces $\mathcal{P}_h$, $\mathcal{W}_h$, $\mathcal{U}_h$, and $\mathcal{Q}_h$, respectively. The basis functions themselves are constructed by
tensor product combinations of nodal and edge polynomials \cite{Gerritsma11, jain2020}. As usual, elements of these spaces can be represented as a linear combination of the basis functions, for example, for $p_{h}\in\mathcal{P}_{h}$, $\boldsymbol{w}_{h}\in\mathcal{W}_{h}$, $\boldsymbol{u}_{h}\in\mathcal{U}_{h}$, and $c_{h}\in\mathcal{Q}_{h}$ we have that
\begin{equation}
	p_{h} = \sum_{i = 1}^{\mathrm{d}_{\mathcal{P}}} \hat{p}_{i}\psi_{i} \,,\quad
\boldsymbol{w}_{h} = \sum_{i = 1}^{\mathrm{d}_{\mathcal{W}}} \hat{w}_{i}\boldsymbol{\alpha}_{i}\,,\quad
\boldsymbol{u}_{h} = \sum_{i = 1}^{\mathrm{d}_{\mathcal{U}}} \hat{u}_{i}\boldsymbol{\beta}_{i}\,, \quad\text{and}\quad
c_{h} = \sum_{i = 1}^{\mathrm{d}_{\mathcal{Q}}} \hat{c}_{i}\gamma_{i}\,. \label{eq:coefficients_1}
\end{equation}
We now introduce an equivalent, but more compact, notation for the expansion of a function as a linear combination of the basis functions
\begin{equation}
	p_{h} = \psi_{h}\valgebra{p}_{h} \,,\quad
\boldsymbol{w}_{h} = \boldsymbol{\alpha}_{h}\valgebra{\boldsymbol{w}}_{h}\,,\quad
\boldsymbol{u}_{h} = \boldsymbol{\beta}_{h}\valgebra{\boldsymbol{u}}_{h}\,, \quad\text{and}
\quad c_{h} = \gamma_{h}\valgebra{c}_{h}\,, \label{eq:coefficients_2}
\end{equation}
where \green{$\valgebra{p}_{h} := [\valgebra{p}_{1}, \dots, \valgebra{p}_{\mathrm{d}_{\mathcal{P}}}]^{\top}$, $\valgebra{\boldsymbol{w}}_{h} := [\valgebra{w}_{1}, \dots, \valgebra{w}_{\mathrm{d}_{\mathcal{W}}}]^{\top}$, $\valgebra{\boldsymbol{u}}_{h} := [\valgebra{u}_{1}, \dots, \valgebra{u}_{\mathrm{d}_{\mathcal{U}}}]^{\top}$, and $\valgebra{c}_{h} := [\valgebra{c}_{1}, \dots, \valgebra{c}_{\mathrm{d}_{\mathcal{Q}}}]^{\top}$} are the column vectors containing the coefficients of the discrete representation of the associated functions, and $\psi_{h} :=  [\psi_{1}, \dots, \psi_{\mathrm{d}_{\mathcal{P}}}]$, $\boldsymbol{\alpha}_{h}:= [\boldsymbol{\alpha}_{1}, \dots, \boldsymbol{\alpha}_{\mathrm{d}_{\mathcal{W}}}]$, $\boldsymbol{\beta}_{h} := [\boldsymbol{\beta}_{1}, \dots, \boldsymbol{\beta}_{\mathrm{d}_{\mathcal{U}}}]$, and $\gamma_{h} := [\gamma_{1}, \dots, \gamma_{\mathrm{d}_{\mathcal{Q}}}]$ are the row vectors containing the basis functions.

\begin{remark}
	Throughout this work we will extensively use the basis functions $\psi_{h}$, $\boldsymbol{\alpha}_{h}$, $\boldsymbol{\beta}_{h}$, and $\gamma_{h}$. For clarity in the notation, these symbols will be solely used to refer to the basis functions:
	\begin{equation}
		\begin{aligned}
			\psi_{h} :=  [\psi_{1}, \dots, \psi_{\mathrm{d}_{\mathcal{P}}}] \qquad\mathrm{such\ that}&\qquad\mathrm{span} \{\psi_1, \dots, \psi_{\mathrm{d}_{\mathcal{P}}}\} = \mathcal{P}_h\\
			\boldsymbol{\alpha}_{h}:= [\boldsymbol{\alpha}_{1}, \dots, \boldsymbol{\alpha}_{\mathrm{d}_{\mathcal{W}}}] \qquad\mathrm{such\ that}&\qquad\mathrm{span} \{\boldsymbol{\alpha}_1, \dots, \boldsymbol{\alpha}_{\mathrm{d}_{\mathcal{W}}}\} = \mathcal{W}_h\\
			\boldsymbol{\beta}_{h} := [\boldsymbol{\beta}_{1}, \dots, \boldsymbol{\beta}_{\mathrm{d}_{\mathcal{U}}}] \qquad\mathrm{such\ that}&\qquad \mathrm{span} \{\boldsymbol{\beta}_1, \dots, \boldsymbol{\beta}_{\mathrm{d}_{\mathcal{U}}}\} = \mathcal{U}_h\\
			\gamma_{h} := [\gamma_{1}, \dots, \gamma_{\mathrm{d}_{\mathcal{Q}}}]\qquad\mathrm{such\ that}&\qquad \mathrm{span} \{\gamma_1, \dots, \gamma_{\mathrm{d}_{\mathcal{Q}}}\} = \mathcal{Q}_h\,.
		\end{aligned}
	\end{equation}
\end{remark}

These discrete function spaces satisfy a
compatibility relation that is expressed as a discrete de Rham complex of the form
\begin{equation}
\mathbb{R}\longrightarrow\mathcal{P}_h \stackrel{\nabla}{\longrightarrow}
\mathcal{W}_h \stackrel{\nabla\times}{\longrightarrow}
\mathcal{U}_h \stackrel{\nabla\cdot}{\longrightarrow} \mathcal{Q}_h \longrightarrow 0.
\end{equation}
The differential operators $\nabla$, $\nabla\times$, and $\nabla\cdot$ are concisely represented at the discrete level by incidence matrices,
$\boldsymbol{\mathsf{E}}^{1,0}$, $\boldsymbol{\mathsf{E}}^{2,1}$, $\boldsymbol{\mathsf{E}}^{3,2}$, such that for $p_{h}\in\mathcal{P}_{h}$, $\boldsymbol{w}_{h}\in\mathcal{W}_{h}$, and $\boldsymbol{u}_{h}\in\mathcal{U}_{h}$ we have that
\begin{equation}
\nabla p_{h} = \boldsymbol{\alpha}_{h}\boldsymbol{\mathsf{E}}^{1,0} \valgebra{p}_{h}, \quad \nabla\times\boldsymbol{w}_{h} = \boldsymbol{\beta}_{h}\boldsymbol{\mathsf{E}}^{2,1} \valgebra{\boldsymbol{w}}_{h}, \quad\text{and}\quad \nabla\cdot\boldsymbol{u}_{h} =  \gamma_{h}\boldsymbol{\mathsf{E}}^{3,2} \valgebra{\boldsymbol{u}}_{h}\,. \label{eq_incidence_matrices}
\end{equation}
In a similar way, we may apply the differential operators to the basis functions $\psi_{h}$, $\boldsymbol{\alpha}_{h}$, $\boldsymbol{\beta}_{h}$, and $\gamma_{h}$ yielding
	\begin{equation}
		\nabla \psi_{h} = \boldsymbol{\alpha}_{h}\boldsymbol{\mathsf{E}}^{1,0}, \quad \nabla\times\boldsymbol{\alpha}_{h} = \boldsymbol{\beta}_{h}\boldsymbol{\mathsf{E}}^{2,1}, \quad\text{and}\quad \nabla\cdot\boldsymbol{\beta}_{h} =  \gamma_{h}\boldsymbol{\mathsf{E}}^{3,2}\,. \label{eq_incidence_matrices_basis}
	\end{equation}
These incidence matrices satisfy the identities
\begin{equation}
	\boldsymbol{\mathsf{E}}^{2,1}\boldsymbol{\mathsf{E}}^{1,0} = \boldsymbol{0}^{2,0}, \quad\text{and}\quad \boldsymbol{\mathsf{E}}^{3,2}\boldsymbol{\mathsf{E}}^{2,1} = \boldsymbol{0}^{3, 1}\,,
\end{equation}
where $\boldsymbol{0}^{2,0}$ is the $\mathrm{d}_{\mathcal{U}} \times \mathrm{d}_{\mathcal{P}}$ zero matrix and $\boldsymbol{0}^{3,1}$ is the $\mathrm{d}_{\mathcal{Q}} \times \mathrm{d}_{\mathcal{W}}$ zero matrix. Note that these two identities are directly related to the well known vector calculus identities $\nabla\times\nabla = 0$ and $\nabla\cdot\nabla\times = 0$. For a more detailed discussion of the incidence matrices and their properties the reader is directed to \cite{Palha14}.

Additionally, these compatible function spaces satisfy an integration-by-parts property, such that the reverse mappings ($\tilde{\nabla}\cdot$, $\tilde{\nabla}\times$, $\tilde{\nabla}$), satisfy
\begin{equation}
0 \longleftarrow\mathcal{P}_h \stackrel{\tilde{\nabla}\cdot}{\longleftarrow}
\mathcal{W}_h \stackrel{\tilde{\nabla}\times}{\longleftarrow}
\mathcal{U}_h \stackrel{\tilde{\nabla}}{\longleftarrow} \mathcal{Q}_h \longleftarrow \mathbb{R}.
\end{equation}
and are represented by weak form adjoint relations. Assuming periodic or homogeneous boundary conditions
we have that for $\boldsymbol{w}_{h}\in\mathcal{W}_{h}$, $\boldsymbol{u}_{h}\in\mathcal{U}_{h}$, and $c_{h}\in\mathcal{Q}_{h}$
\begin{equation}
	\begin{aligned}
		\langle\boldsymbol{\beta}_{h}, \tilde{\nabla} c_{h}\rangle &= -\langle\nabla\cdot\boldsymbol{\beta}_{h}, c_{h}\rangle, \\
		\langle\boldsymbol{\alpha}_{h}, \tilde{\nabla} \times \boldsymbol{u}_{h}\rangle &= \langle\nabla\times\boldsymbol{\alpha}_{h},  \boldsymbol{u}_{h} \rangle, \\
		\langle \psi_{h}, \tilde{\nabla} \cdot \boldsymbol{w}_{h}\rangle &= -\langle  \nabla\psi_{h}, \boldsymbol{w}_{h} \rangle\,.
	\end{aligned} \label{eq:dual_differential_operators}
\end{equation}

\begin{remark}
	In \eqref{eq:dual_differential_operators} we have used an implicit matrix notation. As an example, consider the expression $\langle\boldsymbol{\beta}_{h}, \boldsymbol{\beta}_{h}\rangle$. This expression corresponds to
	\begin{equation}
		\langle\boldsymbol{\beta}_{h}, \boldsymbol{\beta}_{h}\rangle := \int_{\Omega} \boldsymbol{\beta}_{h}^{\top} \boldsymbol{\beta}_{h}\,\mathrm{d}\Omega\,, \label{eq:matrix_inner_product}
	\end{equation}
	which is a matrix of dimensions $\mathrm{d}_{\mathcal{U}}\times\mathrm{d}_{\mathcal{U}}$ with element $(i,j)$, $\langle\boldsymbol{\beta}_{h}, \boldsymbol{\beta}_{h}\rangle_{i,j}$, given by
	\[
		\langle\boldsymbol{\beta}_{h}, \boldsymbol{\beta}_{h}\rangle_{i,j} = \langle\boldsymbol{\beta}_{i}, \boldsymbol{\beta}_{j}\rangle\,.
	\]
	As a more complex example, consider the expression $\langle  \nabla\psi_{h}, \boldsymbol{w}_{h} \rangle$, with $\boldsymbol{w}_{h}\in\mathcal{W}_{h}$. In a similar way, this expression corresponds to
	\[
		\langle  \nabla\psi_{h}, \boldsymbol{w}_{h} \rangle = \int_{\Omega} \nabla\psi_{h}^{\top} \boldsymbol{w}_{h}\,\mathrm{d}\Omega\,,
	\]
	which is a column vector of $\mathrm{d}_{\mathcal{W}}$ elements with element $j$, $\langle  \nabla\psi_{h}, \boldsymbol{w}_{h} \rangle_{j}$, given by \cite{GB16}
	\[
		\langle  \nabla\psi_{h}, \boldsymbol{w}_{h} \rangle_{j} \stackrel{\eqref{eq:matrix_inner_product}}{=} \left(\int_{\Omega}  \nabla\psi_{h}^{\top} \boldsymbol{w}_{h}\right)_{j}\stackrel{\eqref{eq_incidence_matrices_basis}}{=}  \left(\int_{\Omega}  \left(\boldsymbol{\mathsf{E}}^{1,0}\right)^{\top} \boldsymbol{\alpha}^{\top}_{h} \boldsymbol{w}_{h} \right)_{j} =\left(\left(\boldsymbol{\mathsf{E}}^{1,0}\right)^{\top}\langle\boldsymbol{\alpha}_{h}, \boldsymbol{w}_{h} \rangle\right)_{j} = \sum_{k = 1}^{\mathrm{d}_{\mathcal{W}}}\boldsymbol{\mathsf{E}}^{1,0}_{k,j} \langle\boldsymbol{\alpha}_{k}, \boldsymbol{w}_{h}\rangle\,.
	\]
\end{remark}
This leads to the following discrete representations of the dual operators applied to $\boldsymbol{w}_{h}\in\mathcal{P}_{h}$, $\boldsymbol{u}_{h}\in\mathcal{U}_{h}$, and $c_{h}\in\mathcal{Q}_{h}$
\begin{equation}
	\begin{aligned}
		\tilde{\nabla} c_{h} &= -\boldsymbol{\beta}_{h}\,\langle\boldsymbol{\beta}_{h}, \boldsymbol{\beta}_{h}\rangle^{-1} \left(\boldsymbol{\mathsf{E}}^{3,2}\right)^{\top} \langle\gamma_{h}, \gamma_{h}\rangle\,\valgebra{c}_{h}\,, \\
		\tilde{\nabla}\times \boldsymbol{u}_{h} &= \boldsymbol{\alpha}_{h}\,\langle\boldsymbol{\alpha}_{h}, \boldsymbol{\alpha}_{h}\rangle^{-1} \left(\boldsymbol{\mathsf{E}}^{2,1}\right)^{\top} \langle\boldsymbol{\beta}_{h}, \boldsymbol{\beta}_{h}\rangle\,\valgebra{\boldsymbol{u}}_{h}\,,\\
		\tilde{\nabla}\cdot \boldsymbol{w}_{h} &= -\psi_{h}\, \langle\psi_{h}, \psi_{h}\rangle^{-1} \left(\boldsymbol{\mathsf{E}}^{1,0}\right)^{\top} \langle\boldsymbol{\alpha}_{h}, \boldsymbol{\alpha}_{h}\rangle\, \valgebra{\boldsymbol{w}}_{h}\,.
	\end{aligned}
\end{equation}

The discrete form of the variational 3D compressible Euler equations \eqref{eq::ce_var} may then be expressed as: Given $\left.\boldsymbol{u}_{h}\right|_{t=0}\in\mathcal{U}$, $\left.\rho\right|_{t=0}, \left.\Theta\right|_{t=0}\in\mathcal{Q}$, find $\boldsymbol{u}_{h}\in\mathcal{U}$, $\rho, \Theta\in\mathcal{Q}$ such that
\begin{subequations}\label{eq::ce_var_disc}
\begin{align}
\Bigg\langle\boldsymbol{\beta}_{h},\frac{\partial\boldsymbol{u}_{h}}{\partial t}\Bigg\rangle &=
-\langle\boldsymbol{\beta}_{h},\boldsymbol{q}_{h}\times\boldsymbol{U}_{h}\rangle
+\langle\nabla\cdot\boldsymbol{\beta}_{h},\Phi_{h}\rangle
+\langle\nabla\cdot\mathbb{P}^{\mathrm{div}}\left[\theta_{h}\boldsymbol{\beta}_{h}\right], \Pi_{h}\rangle\,,\label{eq::ce_var_mom_disc}\\
\Bigg\langle\gamma_{h},\frac{\partial\rho_{h}}{\partial t}\Bigg\rangle &=
-\langle\gamma_{h},\nabla\cdot\boldsymbol{U}_{h}\rangle\,,\\
\Bigg\langle\gamma_{h},\frac{\partial\Theta_{h}}{\partial t}\Bigg\rangle &=
-\langle\gamma_{h},\nabla\cdot\mathbb{P}^{\mathrm{div}}\left[\theta_{h}\boldsymbol{U}_{h}\right]\rangle\,,
\end{align}
\end{subequations}
and $\boldsymbol{U}_{h} \in \mathcal{U}$, $\Phi_{h}, \Pi_{h} \in \mathcal{Q}$ such that
\begin{subequations} \label{eq::ce_prognostic}
\begin{align}
	\langle\boldsymbol{\beta}_h,\frac{\delta H_h}{\delta\boldsymbol{u}_h}\rangle &= \langle\boldsymbol{\beta}_{h}, \boldsymbol{U}_{h}\rangle = \langle \boldsymbol{\beta}_{h}, \rho_{h}\boldsymbol{u}_{h} \rangle\,, \\
	\langle\gamma_h,\frac{\delta H_h}{\delta\rho_h}\rangle &= \langle\gamma_{h}, \Phi_{h}\rangle = \langle \gamma_{h}, \frac{1}{2}\boldsymbol{u}_{h}\cdot\boldsymbol{u}_{h} + gz\rangle\,, \\
	\langle\gamma_h,\frac{\delta H_h}{\delta\Theta_h}\rangle &= \langle\gamma_{h}, \Pi_{h}\rangle = \langle \gamma_{h}, c_p\Bigg(\frac{R\Theta_h}{p_0}\Bigg)^{\frac{R}{c_v}}\rangle.
\end{align}
\end{subequations}

The system of equations \eqref{eq::ce_var_disc}, obtained after spatial discretisation, constitutes a system of ordinary differential equations. 
%The objective we propose to achieve now it 
\green{We now wish }
to devise a discrete time integration formulation that preserves the energy conservation properties of 
\green{the spatial discretisation in }
\eqref{eq::ce_var_disc}. To introduce the energy conserving time integration scheme employed in this work, 
we recall the total energy $H$, \eqref{eq::total_energy},
\begin{equation}
{H}_{h} = \int_{\Omega}\overbrace{\frac{1}{2}\rho_{h}\boldsymbol{u}_{h}\cdot\boldsymbol{u}_{h}}^{{K_{h}}} + \overbrace{\vphantom{\frac{1}{2}}\rho_{h} gz}^{{P_{h}}} + \overbrace{\frac{c_v}{c_p}\Theta_{h}\Pi_{h}}^{{I_{h}}}\mathrm{d}\Omega,
\end{equation}
and its time derivative, \eqref{eq::hamiltonian_time_variation_0},
\begin{equation}\label{eq::discrete_hamiltonian}
\Bigg\langle\frac{\delta H_{h}}{\delta\boldsymbol{u}_{h}}, \frac{\partial\boldsymbol{u}_{h}}{\partial t}\Bigg\rangle +
\Bigg\langle\frac{\delta H_{h}}{\delta\rho_{h}},\frac{\partial\rho_{h}}{\partial t}\Bigg\rangle +
\Bigg\langle\frac{\delta H_{h}}{\delta\Theta_{h}},\frac{\partial\Theta_{h}}{\partial t}\Bigg\rangle =
\Bigg\langle\boldsymbol{U}_{h}, \frac{\partial\boldsymbol{u}_{h}}{\partial t}\Bigg\rangle +
\Bigg\langle\Phi_{h},\frac{\partial\rho_{h}}{\partial t}\Bigg\rangle +
\Bigg\langle\Pi_{h},\frac{\partial\Theta_{h}}{\partial t}\Bigg\rangle =
\frac{\mathrm{d}H_{h}}{\mathrm{d} t} = 0.
\end{equation}
We now note that \eqref{eq::ce_var_disc} can be written in matrix form as
\begin{equation}
	\def\arraystretch{2.2}
	\left[
	\begin{array}{>{\displaystyle}c}
		\langle\boldsymbol{\beta}_{h}, \boldsymbol{\beta}_{h}\rangle\frac{\mathrm{d}\valgebra{\boldsymbol{u}}_{h}}{\mathrm{d} t} \\
		\langle\gamma_{h}, \gamma_{h}\rangle\frac{\mathrm{d}\valgebra{\rho}_{h}}{\mathrm{d} t} \\
		\langle\gamma_{h}, \gamma_{h}\rangle\frac{\mathrm{d}\valgebra{\Theta}_{h}}{\mathrm{d} t}
	\end{array}
	\right]
	=
	\left[
		\begin{array}{ccc}
			-\langle\boldsymbol{\beta}_{h},\boldsymbol{q}_{h}\times\boldsymbol{\beta}_{h}\rangle & \langle\nabla\cdot\boldsymbol{\beta}_{h},\gamma_{h}\rangle & \langle\nabla\cdot\mathbb{P}^{\mathrm{div}}\left[\theta_{h}\boldsymbol{\beta}_{h}\right], \gamma_{h}\rangle \\
			-\langle\gamma_{h},\nabla\cdot\boldsymbol{\beta}_{h}\rangle & 0 & 0 \\
			-\langle\gamma_{h},\nabla\cdot\mathbb{P}^{\mathrm{div}}\left[\theta_{h}\boldsymbol{\beta}_{h}\right]\rangle & 0 & 0
		\end{array}
	\right]
	\left[
	\begin{array}{c}
		\hat{\boldsymbol{U}}_{h} \\
		\hat{\Phi}_{h} \\
		\hat{\Pi}_{h}
	\end{array}
	\right]\,. \label{eq:skew_symmetric}
\end{equation}

\subsection{Temporal discretisation}

Equation \eqref{eq:skew_symmetric} is essentially a system of ordinary differential equations (ODEs) with the general form
\begin{equation}
	\boldsymbol{\mathsf{M}}\frac{\mathrm{d}\valgebra{\boldsymbol{y}}(t)}{\mathrm{d} t} = \boldsymbol{\mathsf{S}}(\valgebra{\boldsymbol{y}}(t))\, \valgebra{\boldsymbol{x}}(t), \label{eq:skew_symmetric_simplified}
\end{equation}
where
\begin{equation}
	\valgebra{\boldsymbol{y}}(t) := 
	\left[
	\begin{array}{c}
		\valgebra{\boldsymbol{u}}_{h} \\
		\valgebra{\rho}_{h} \\
		\valgebra{\Theta}_{h}
	\end{array}
	\right]\,,\qquad  \valgebra{\boldsymbol{x}}(t) :=  
	\left[
	\begin{array}{c}
		\hat{\boldsymbol{U}}_{h} \\
		\hat{\Phi}_{h} \\
		\hat{\Pi}_{h}
	\end{array}
	\right]\,, \label{eq:y_definition}
\end{equation}
\green{$\boldsymbol{\mathsf{S}}(\valgebra{\boldsymbol{y}}(t))$ is a time dependent skew-symmetric matrix
%\green{the time dependent skew-symmetric matrix $\boldsymbol{\mathsf{S}}(\valgebra{\boldsymbol{y}}(t))$ is
%\begin{equation}
	%\boldsymbol{\mathsf{S}}(\valgebra{\boldsymbol{y}}(t)) := \left[
		%\begin{array}{ccc}
			%-\langle\boldsymbol{\beta}_{h},\boldsymbol{q}_{h}\times\boldsymbol{\beta}_{h}\rangle & \langle\nabla\cdot\boldsymbol{\beta}_{h},\gamma_{h}\rangle & \langle\nabla\cdot\mathbb{P}^{\mathrm{div}}\left[\theta_{h}\boldsymbol{\beta}_{h}\right], \gamma_{h}\rangle \\
			%-\langle\gamma_{h},\nabla\cdot\boldsymbol{\beta}_{h}\rangle & 0 & 0 \\
			%-\langle\gamma_{h},\nabla\cdot\mathbb{P}^{\mathrm{div}}\left[\theta_{h}\boldsymbol{\beta}_{h}\right]\rangle & 0 & 0
		%\end{array}
	%\right]\,,
%\end{equation}
%and the constant in time symmetric positive definite matrix $\boldsymbol{\mathsf{M}}$ is
and $\boldsymbol{\mathsf{M}}$ is a constant in time symmetric positive definite matrix.
%\begin{equation}
	%\boldsymbol{\mathsf{M}} := \left[
		%\begin{array}{ccc}
			%\langle\boldsymbol{\beta}_{h}, \boldsymbol{\beta}_{h}\rangle & 0 & 0 \\
			%0 & \langle\gamma_{h}, \gamma_{h}\rangle & 0 \\
			%0 & 0 & \langle\gamma_{h}, \gamma_{h}\rangle
		%\end{array}
	%\right]\,.
%\end{equation}
%The variable $\hat{\boldsymbol{y}} := [\hat{\boldsymbol{u}}_{h}, \hat{\rho}_{h}, \hat{\Theta}_{h}]^{\top}$, 
The variable $\hat{\boldsymbol{y}}$  
\eqref{eq:y_definition} is a column vector containing the coefficients of the discrete representation of the fields $\boldsymbol{u}_{h}$, $\rho_{h}$, and $\Theta_{h}$, as introduced in \eqref{eq:coefficients_1} and \eqref{eq:coefficients_2}. In this way we have that
\begin{equation}
	\hat{y}_{i} = \hat{u}_{i}, \quad \mathrm{for}\quad i = 1, \dots, \mathrm{d}_{\mathcal{U}}, \quad \hat{y}_{i + \mathrm{d}_{\mathcal{U}}} = \hat{\rho}_{i}, \quad \mathrm{for}\quad i = 1, \dots, \mathrm{d}_{\mathcal{Q}}, \quad\mathrm{and}\quad \hat{y}_{i + \mathrm{d}_{\mathcal{U}} + \mathrm{d}_{\mathcal{Q}}} = \hat{\Theta}_{i}, \quad \mathrm{for}\quad i = 1, \dots, \mathrm{d}_{\mathcal{Q}}\,.
\end{equation} 
These coefficients $\hat{\boldsymbol{y}}$ are time dependent and their evolution is governed by the system of ODEs \eqref{eq:skew_symmetric_simplified}.
 
To solve this system of ODEs in time, we follow a similar approach to that for the spatial discretisation. 
We approximate the time dependent coefficients $\valgebra{\boldsymbol{y}}(t)$ on the time step interval $[t^{n}, t^{n+1}]$ with a temporal discrete polynomial function space of degree $s$, $\mathcal{T}^{0}_{h}([t^{n}, t^{n+1}]) := \mathrm{span}\{l_{0}, \dots, l_{s}\}$ where $l_{i}(t)$, $i=0,\dots s$, are the ($s+1$) nodal polynomials of degree $s$, as mentioned before for the spatial discretisation. In this way the polynomial temporal expansion of the spatial coefficients $\valgebra{y}_{i}^{h}$ is
\begin{equation}
	\valgebra{y}_{i}^{h} = l_{h}\, \valgebra{\valgebra{y}}_{i}^{h} := \sum_{j=0}^{s} l_{j}(t) \valgebra{\valgebra{y}}^{j}_{i}\,, \qquad t\in([t^{n}, t^{n+1}])\,,\label{eq:temporal_notation}
\end{equation}
where, as in \eqref{eq:coefficients_1} and \eqref{eq:coefficients_2}, 
$\valgebra{\valgebra{y}}_{i}^{h} := [\valgebra{\valgebra{y}}_{i}^{1}, \dots \valgebra{\valgebra{y}}_{i}^{s+1}]$. The double hat notation 
$\valgebra{\valgebra{y}}_{i}^{h}$ in \eqref{eq:temporal_notation} highlights the double discretization taking place (spatial and temporal). 
The terms $\valgebra{y}_{i}^{h}$ with $i = 1, \dots, \mathrm{d}_{\mathcal{U}} + 2\mathrm{d}_{\mathcal{Q}}$ correspond to the time dependent 
polynomial expansion of the $i$-th coefficient of the spatial polynomial expansion $\hat{\boldsymbol{y}} := [\hat{\boldsymbol{u}}_{h}, \hat{\rho}_{h}, \hat{\Theta}_{h}]^{\top}$. The terms $\valgebra{\valgebra{y}}^{j}_{i}$ correspond to the $j$-th coefficient of the temporal polynomial expansion associated to the $i$-th coefficient of the spatial polynomial expansion.}

The time dependent coefficients $\valgebra{\boldsymbol{x}}(t)$ are approximated with a temporal discrete polynomial function space of degree ($s-1$), $\mathcal{T}^{1}_{h}([t^{n}, t^{n+1}]) := \mathrm{span}\{e_{1}, \dots, e_{s}\}$ where $e_{i}(t)$, $i=1,\dots s$, are the $s$ edge polynomials of degree ($s-1$), as mentioned before for the spatial discretisation. In this way
\green{\begin{equation}
	\valgebra{x}_{i}^{h} = e_{h}\, \valgebra{\valgebra{x}}_{i}^{h} := \sum_{j=1}^{s} e_{j}(t) \valgebra{\valgebra{x}}^{j}_{i}\,, \qquad t\in([t^{n}, t^{n+1}])\,.
\end{equation}}
As before, these two function spaces constitute a discrete de Rham complex
\begin{equation}\label{eq:de_rham_time}
\mathbb{R}\longrightarrow\mathcal{T}^{0}_{h}([t^{n}, t^{n+1}]) \stackrel{\frac{\mathrm{d}}{\mathrm{d}t}}{\longrightarrow}
\mathcal{T}^{1}_{h}([t^{n}, t^{n+1}]) \longrightarrow 0\,,
\end{equation}
and we have that
\green{\begin{equation}
	\frac{\mathrm{d}\valgebra{y}_{i}^{h}}{\mathrm{d}t}  = e_{h}\boldsymbol{\mathsf{E}}^{1,0}_{t}\valgebra{\valgebra{y}}_{i}^{h}\,,\qquad\mathrm{and}\qquad \frac{\mathrm{d}l_{h}}{\mathrm{d}t} =  e_{h}\boldsymbol{\mathsf{E}}^{1,0}_{t}\,.
\end{equation}}
\begin{remark}
	Note that here we added the subscript $t$ to the incidence matrix, $\boldsymbol{\mathsf{E}}^{1,0}_{t}$, to distinguish it from the spatial ones. \green{Also, as mentioned before, the double hat highlights the fact that we have performed both discretisation in space and time and these coefficients are associated to this double discretisation. } When an explicit reference to the temporal indices is required they are presented as superscript, as in \eqref{eq:temporal_notation}.
\end{remark}

\begin{remark}
	The discrete de Rham complex \green{\eqref{eq:de_rham_time} } is associated to a single time step (for simplicity of exposition and computational efficiency). Naturally, as for the spatial discretisation, considering multiple time steps is a straightforward extension. In this work we perform the discretisation in a per step fashion.	
\end{remark}

The temporal discretisation then follows a similar mixed Galerkin approach as presented for the spatial discretisation: Given \green{$\valgebra{x}_{i}^{h} \in \mathcal{T}^{1}_{h}([t^{n}, t^{n+1}])$, with $i = 1, \dots, (\mathrm{d}_{\mathcal{U}} + 2\mathrm{d}_{\mathcal{Q}})$}, find \green{$\valgebra{y}_{i}^{h} \in \mathcal{T}^{0}_{h}([t^{n}, t^{n+1}])$, with $i = 1, \dots, (\mathrm{d}_{\mathcal{U}} + 2\mathrm{d}_{\mathcal{Q}})$}, such that
\begin{equation}
	\boldsymbol{\mathsf{M}}\,\langle e_{h}, \frac{\mathrm{d}\valgebra{\boldsymbol{y}}_{h}(t)}{\mathrm{d} t}\rangle_{[t^{n}, t^{n+1}]} = \langle e_{h}, \boldsymbol{\mathsf{S}}(\valgebra{\boldsymbol{y}}_{h}(t))\, \valgebra{\boldsymbol{x}}_{h}(t)\rangle_{[t^{n}, t^{n+1}]}\,. \label{eq:innerproduct_time}
\end{equation}

\begin{remark}
	The time-related inner products follow the same rules outlined before, namely
	\begin{equation}
		\langle e_{h}, \frac{\mathrm{d}\valgebra{\boldsymbol{y}}_{h}(t)}{\mathrm{d} t}\rangle_{[t^{n}, t^{n+1}]} := \int_{t^{n}}^{t^{n+1}} e_{h}^{\top} \frac{\mathrm{d}\valgebra{\boldsymbol{y}}_{h}(t)}{\mathrm{d} t} \,\mathrm{d}t\,.
	\end{equation}
	Therefore, \eqref{eq:innerproduct_time} is
	\begin{equation}
		\boldsymbol{\mathsf{M}}\,\int_{t^{n}}^{t^{n+1}} e_{h}^{\top} \frac{\mathrm{d}\valgebra{\boldsymbol{y}}_{h}(t)}{\mathrm{d} t}\,\mathrm{d}t = \int_{t^{n}}^{t^{n+1}} e_{h}^{\top} \boldsymbol{\mathsf{S}}(\valgebra{\boldsymbol{y}}_{h}(t))\, \valgebra{\boldsymbol{x}}_{h}(t)\,\mathrm{d}t\,. \label{eq:long_time_integrator}
	\end{equation}
Note that we implicitly separate the spatial from the temporal coefficients, for simplicity of the notation, i.e., the constant in time 
matrix $\boldsymbol{\mathsf{M}}$ and the time dependent matrix  $\boldsymbol{\mathsf{S}}(\valgebra{\boldsymbol{y}}_{h}(t))$ apply at the 
spatial level also.
\end{remark}

This choice of time integration is notable because it is energy conserving. To prove this, we recall the expression for the time variation 
of the energy of this system of equations (the Hamiltonian, $H_{h}$) \eqref{eq::discrete_hamiltonian} 
\begin{equation}
\int_{t^{n}}^{t^{n+1}}\frac{\mathrm{d}H_{h}}{\mathrm{d} t}\,\mathrm{d}t = H_{h}(\boldsymbol{y}_{h}(t^{n+1})) - H_{h}(\boldsymbol{y}_{h}(t^{n})) = \int_{t^{n}}^{t^{n+1}}\Bigg\langle\frac{\delta{H_{h}}}{\delta\boldsymbol{y}_{h}},\frac{\partial\boldsymbol{y}_{h}}{\partial t}\Bigg\rangle\,\mathrm{d}t. \label{eq::hamiltonian_time_variation}
\end{equation}
We may now use \eqref{eq::ce_prognostic} to rewrite \eqref{eq::hamiltonian_time_variation} as
\begin{equation}
\int_{t^{n}}^{t^{n+1}}\frac{\mathrm{d}H_{h}}{\mathrm{d} t}\,\mathrm{d}t = 
H_{h}(\boldsymbol{y}_{h}(t^{n+1})) - H_{h}(\boldsymbol{y}_{h}(t^{n})) = 
\int_{t^{n}}^{t^{n+1}}\Bigg\langle\boldsymbol{x}_{h},\frac{\partial\boldsymbol{y}_{h}}{\partial t}\Bigg\rangle\,\mathrm{d}t = 
\int_{t^{n}}^{t^{n+1}}\valgebra{\boldsymbol{x}}_{h}^{\top}\boldsymbol{\mathsf{M}}\frac{\mathrm{d}\valgebra{\boldsymbol{y}}_{h}}{\mathrm{d} t}\,\mathrm{d}t = 
\valgebra{\valgebra{\boldsymbol{x}}}_{h}^{\top}\boldsymbol{\mathsf{M}}\int_{t^{n}}^{t^{n+1}}e_{h}^{\top}\frac{\mathrm{d}\valgebra{\boldsymbol{y}}_{h}}{\mathrm{d} t}\,\mathrm{d}t  . \label{eq::hamiltonian_time_variation_2}
\end{equation}
Finally, using the proposed time integration scheme, \eqref{eq:long_time_integrator}, we obtain energy conservation
\begin{equation}
H_{h}(\boldsymbol{y}_{h}(t^{n+1})) - H_{h}(\boldsymbol{y}_{h}(t^{n})) = 
\valgebra{\valgebra{\boldsymbol{x}}}_{h}^{\top}\boldsymbol{\mathsf{M}}\int_{t^{n}}^{t^{n+1}}e_{h}^{\top}
e_{h}\,\mathrm{d}t(\valgebra{\valgebra{\boldsymbol{y}}}_{h}^{n+1} - \valgebra{\valgebra{\boldsymbol{y}}}_{h}^{n}) = 
\valgebra{\valgebra{\boldsymbol{x}}}_{h}^{\top}\int_{t^{n}}^{t^{n+1}} e_{h}^{\top} \boldsymbol{\mathsf{S}}(\valgebra{\boldsymbol{y}}_{h}(t))\, \valgebra{\boldsymbol{x}}_{h}(t)\,\mathrm{d}t = 
\valgebra{\valgebra{\boldsymbol{x}}}_{h}^{\top}\left(\int_{t^{n}}^{t^{n+1}} e_{h}^{\top} \boldsymbol{\mathsf{S}}(\valgebra{\boldsymbol{y}}_{h}(t))\,  e_{h}\,\mathrm{d}t\right) \, \valgebra{\valgebra{\boldsymbol{x}}}_{h}. \label{eq::hamiltonian_time_variation_3}
\end{equation}
Since the term inside the parenthesis on the right hand side is skew-symmetric, the quadratic term cancels and we obtain energy conservation at the discrete level
\begin{equation}
H_{h}(\boldsymbol{y}_{h}(t^{n+1})) - H_{h}(\boldsymbol{y}_{h}(t^{n})) = 0. \label{eq::hamiltonian_time_variation_4}
\end{equation}

\begin{remark}
	 In equation \eqref{eq::hamiltonian_time_variation_3}, the term $\int_{t^{n}}^{t^{n+1}}e_{h}^{\top} e_{h}\,\mathrm{d}t$ is a temporal mass matrix.
\end{remark}

\begin{remark}
	Two key ingredients need to be highlighted. The first is the construction of a system of equations of the form \eqref{eq:skew_symmetric} and \eqref{eq:skew_symmetric_simplified} where $\boldsymbol{\mathsf{S}}$ is skew-symmetric. This was only possible to obtain by a judicious choice of discrete function spaces and an equally careful selection of physical field quantities to employ in the construction of the system of equations. The second key ingredient is the construction of the time integration scheme also based on a sequence of polynomial spaces together with a Galerkin projection employing exact integration. Exact temporal integration is fundamental to guarantee the equality needed to perform the last step, equation \eqref{eq::hamiltonian_time_variation_3}.
\end{remark}
This temporal discretisation is very closely related to the works by Hairer et al. \cite{Hairer2010, Hairer2011}. For the lowest order case (discussed below) this approach results in time stepping scheme identical to \cite{Hairer2011}.

In this formulation any polynomial degree $s$ may be employed to approximate the solution in time. In this work we choose $s = 1$ since this choice will greatly simplify the temporal discretisation. If this choice is made we will have
\begin{equation}
	l_{0}(t) = \frac{t^{n+1} - t}{\Delta t}\,, \quad l_{1}(t) = \frac{t - t^{n}}{\Delta t}\,, \quad \mathrm{and}\quad e_{1}(t) = \frac{1}{\Delta t}\,.
\end{equation}
where $\Delta t := t^{n+1} - t^{n}$. \green{In this case, since $e_{h} = e_{1}(t) = \frac{1}{\Delta t}$, the coefficients of the temporal polynomial expansion of $\valgebra{\boldsymbol{x}}_{h}(t)$ become
\begin{equation}
	\valgebra{\valgebra{\boldsymbol{x}}}_{h}(t) = [\valgebra{\valgebra{U}}^{1}_{1}, \dots, \valgebra{\valgebra{U}}^{1}_{\mathrm{d}_{\mathcal{U}}}, \valgebra{\valgebra{\Phi}}^{1}_{1}, \dots, \valgebra{\valgebra{\Phi}}^{1}_{\mathrm{d}_{\mathcal{Q}}}, \valgebra{\valgebra{\Pi}}^{1}_{1}, \dots, \valgebra{\valgebra{\Pi}}^{1}_{\mathrm{d}_{\mathcal{Q}}}]^{\top}\,e_{1}(t)\,. \label{eq::x_coefficients}
\end{equation}
We recall now the properties of the edge polynomials, \cite{Gerritsma11, jain2020}, and its associated coefficients. 
Just as the nodal polynomial expansion interpolates a function, the edge polynomial expansion histopolates (interpolates 
the integral of) a function. While the coefficients of the nodal polynomial expansion correspond to pointwise sampling of a 
function, the coefficients of the edge polynomial expansion correspond to sampling of the integral of the function over 
intervals \cite{Gerritsma11, jain2020}. Note that in the temporal discretization considered here there is only one polynomial 
basis in the edge polynomial expansion and therefore the temporal expansion of each spatial coefficient has only one coefficient. 
Given the integral interpolation properties of the edge basis functions, the single coefficient of this expansion is the temporal integral
\begin{equation}
	\valgebra{\valgebra{\boldsymbol{x}}}_{i}^{1} = \valgebra{\valgebra{U}}^{1}_{i} = \int_{t^{n}}^{t^{n+1}}\valgebra{U}_{i}^{h}(t) \,\mathrm{d}t, \quad \valgebra{\valgebra{\boldsymbol{x}}}_{j + \mathrm{d}_{\mathcal{U}}}^{1} = \valgebra{\valgebra{\Phi}}^{1}_{j} = \int_{t^{n}}^{t^{n+1}}\valgebra{\Phi}_{j}^{h}(t) \,\mathrm{d}t, \quad \valgebra{\valgebra{\boldsymbol{x}}}_{k + \mathrm{d}_{\mathcal{U}} + \mathrm{d}_{\mathcal{Q}}}^{1} = \valgebra{\valgebra{\Pi}}^{1}_{k} = \int_{t^{n}}^{t^{n+1}}\valgebra{\Pi}_{k}^{h}(t) \,\mathrm{d}t\,, \label{eq:mean_01}
\end{equation}
with $i = 1, \dots, \mathrm{d}_{\mathcal{U}}$, and $j, k = 1, \dots, \mathrm{d}_{\mathcal{Q}}$. We now introduce 
\begin{equation}
	\overline{\valgebra{\valgebra{\boldsymbol{x}}}_{i}^{h}} := \frac{\valgebra{\valgebra{\boldsymbol{x}}}_{i}^{1}}{\Delta t} = \frac{\int_{t^{n}}^{t^{n+1}}\valgebra{\valgebra{\boldsymbol{x}}}_{i}^{h}\, \mathrm{d}t}{\Delta t}\,,\label{eq:mean_02}
\end{equation}
the exact time integral of the temporal approximation between time levels
$n$ and $n+1$. Following the notation introduced before, \eqref{eq:coefficients_1}, we will use 
\begin{equation}
	\overline{\valgebra{\valgebra{\boldsymbol{x}}}}_{h} := [\overline{\valgebra{\valgebra{\boldsymbol{x}}}_{1}^{h}}, \dots, \overline{\valgebra{\valgebra{\boldsymbol{x}}}^{h}}_{\mathrm{d}_{\mathcal{U}} + 2\mathrm{d}_{\mathcal{Q}}}]^{\top} \stackrel{\eqref{eq::x_coefficients} + \eqref{eq:mean_01} + \eqref{eq:mean_02}}{=} \frac{1}{\Delta t}[\valgebra{\valgebra{U}}^{1}_{1}, \dots, \valgebra{\valgebra{U}}^{1}_{\mathrm{d}_{\mathcal{U}}}, \valgebra{\valgebra{\Phi}}^{1}_{1}, \dots, \valgebra{\valgebra{\Phi}}^{1}_{\mathrm{d}_{\mathcal{Q}}}, \valgebra{\valgebra{\Phi}}^{1}_{1}, \dots, \valgebra{\valgebra{\Pi}}^{1}_{\mathrm{d}_{\mathcal{Q}}}]^{\top}\,.
\end{equation}}
If we substitute this into \eqref{eq:long_time_integrator} we obtain
\green{\begin{equation}
		\boldsymbol{\mathsf{M}}\,\left(\valgebra{\valgebra{\boldsymbol{y}}}^{n+1}_{h} - \valgebra{\valgebra{\boldsymbol{y}}}^{n}_{h}\right) = \left(\int_{t^{n}}^{t^{n+1}} \boldsymbol{\mathsf{S}}(\valgebra{\boldsymbol{y}}_{h}(t))\,\mathrm{d}t\right) \, \frac{\valgebra{\valgebra{\boldsymbol{x}}}_{h}}{\Delta t} \stackrel{\eqref{eq:mean_02}}{=} \left(\int_{t^{n}}^{t^{n+1}} \boldsymbol{\mathsf{S}}(\valgebra{\boldsymbol{y}}_{h}(t))\,\mathrm{d}t\right) \, \overline{\valgebra{\valgebra{\boldsymbol{x}}}}_{h}\,. \label{eq:long_time_integrator_s_1}
	\end{equation}}
	Since $\boldsymbol{\mathsf{S}}(\valgebra{\boldsymbol{y}}_{h}(t))$ depends only linearly on time, we may use lowest order Gauss integration (midpoint rule) to exactly integrate the term involving $\boldsymbol{\mathsf{S}}(\valgebra{\boldsymbol{y}}_{h}(t))$, yielding
	\begin{equation}
		\boldsymbol{\mathsf{M}}\,\left(\valgebra{\valgebra{\boldsymbol{y}}}^{n+1}_{h} - \valgebra{\valgebra{\boldsymbol{y}}}^{n}_{h}\right) = \Delta t \,\boldsymbol{\mathsf{S}}\left(\frac{\valgebra{\valgebra{\boldsymbol{y}}}^{n+1}_{h} + \valgebra{\valgebra{\boldsymbol{y}}}^{n}_{h}}{2}\right) \, \overline{\valgebra{\valgebra{\boldsymbol{x}}}}_{h}\,. \label{eq:long_time_integrator_s_2}
	\end{equation}
	
Therefore this system of equations may be cast in a fully discrete form as:
\begin{multline}
	\left[
		\begin{array}{c}
			\langle\boldsymbol{\beta}_{h}, \boldsymbol{\beta}_{h}\rangle\,\hat{\boldsymbol{u}}_{h}^{n+1} \\
			\langle\gamma_{h}, \gamma_{h}\rangle\,\hat{\rho}_{h}^{n+1} \\
			\langle\gamma_{h}, \gamma_{h}\rangle\,\hat{\Theta}_{h}^{n+1}
		\end{array}
	\right]
	=
	\left[
		\begin{array}{c}
			\langle\boldsymbol{\beta}_{h}, \boldsymbol{\beta}_{h}\rangle\,\hat{\boldsymbol{u}}_{h}^{n} \\
			\langle\gamma_{h}, \gamma_{h}\rangle\,\hat{\rho}_{h}^{n} \\
			\langle\gamma_{h}, \gamma_{h}\rangle\,\hat{\Theta}_{h}^{n}
		\end{array}
	\right]
	- \\
	\Delta t
	\left[
		\begin{array}{ccc}
			\langle\boldsymbol{\beta}_{h},\boldsymbol{q}^{n+\frac{1}{2}}_{h}\times\boldsymbol{\beta}_{h}\rangle  & -\left(\boldsymbol{\mathsf{E}}^{3,2}\right)^{\top}\langle\gamma_{h}, \gamma_{h}\rangle & -\langle\boldsymbol{\beta}_{h}, \theta^{n+\frac{1}{2}}_{h}\boldsymbol{\beta}_{h}\rangle \langle\boldsymbol{\beta}_{h}, \boldsymbol{\beta}_{h}\rangle^{-1}\left(\boldsymbol{\mathsf{E}}^{3,2}\right) \langle\gamma_{h}, \gamma_{h}\rangle \\
			\langle\gamma_{h}, \gamma_{h}\rangle\boldsymbol{\mathsf{E}}^{3,2} & 0 & 0 \\
			\langle\gamma_{h}, \gamma_{h}\rangle\boldsymbol{\mathsf{E}}^{3,2} \langle\boldsymbol{\beta}_{h}, \boldsymbol{\beta}_{h}\rangle^{-1}\langle\boldsymbol{\beta}_{h}, \theta^{n + \frac{1}{2}}_{h}\boldsymbol{\beta}_{h}\rangle   & 0 & 0
		\end{array}
	\right]
	\left[
		\begin{array}{c}
			\overline{\valgebra{\boldsymbol{U}}}_{h} \\
			\overline{\valgebra{\Phi}}_{h} \\
			\overline{\valgebra{\Pi}}_{h}
		\end{array}
	\right]\,.  \label{eq::de}
\end{multline}
\green{As mentioned before,  $\overline{{\valgebra{\boldsymbol{U}}}}_h$, $\overline{\valgebra{{\Phi}}}_h$,
and $\overline{\valgebra{{\Pi}}}_h$ are the exact time integrals of the temporal polynomial reconstruction of $\valgebra{\boldsymbol{U}}_{h}$, $\valgebra{\Phi}_{h}$, and $\valgebra{\Pi}_{h}$. An important point to note is the following. Since $\valgebra{U}_{h}$, $\valgebra{\Phi}_{h}$, and $\valgebra{\Pi}_{h}$ are functions of $\valgebra{\boldsymbol{u}}$, $\valgebra{\rho}$, and $\valgebra{\Theta}$, \eqref{eq::variational_derivatives_hamiltonian}, we employ
discrete piecewise linear approximations between time levels
$n$ and $n+1$ to be described in the following section,
and ${\theta}^{n + \frac{1}{2}}_{h} = (\theta_h^n + \theta_h^{n+1})/2$,
${\boldsymbol{q}}^{n + \frac{1}{2}}_{h} = (\boldsymbol{q}_h^n + \boldsymbol{q}_h^{n+1})/2$
are the time centered potential temperature and potential vorticity.}

\section{HEVI splitting}\label{sec::hevi}

The above system \eqref{eq::de} may be dimensionally split into an implicit solve for the vertical dynamics (incorporating the
horizontal divergence terms) and explicit momentum advection by defining the horizontal and vertical velocity components
respectively as $\boldsymbol{v}_h\in\mathcal{U}_h^{\parallel}$, $\boldsymbol{w}_h\in\mathcal{U}_h^{\perp}$, such that
$\boldsymbol{u}_h = \boldsymbol{v}_h + \boldsymbol{w}_h\in\mathcal{U}_h$ \cite{LP20}. Similarly \green{the discrete vorticity
vector, $\boldsymbol{q}_h\in\mathcal{W}_h$ may be partitioned into its respective vertical and horizontal components (in global 
coordinates), $q_h^{\perp} = \boldsymbol{q}_h\cdot\boldsymbol{e}_z$, 
$\boldsymbol{q}_h^{\parallel} = (-\boldsymbol{q}_h\cdot\boldsymbol{e}_{\phi},\boldsymbol{q}\cdot\boldsymbol{e}_{\lambda})$, where 
$\boldsymbol{e}_{\lambda}$, $\boldsymbol{e}_{\phi}$, $\boldsymbol{e}_z$ are unit vectors in the zonal, meridional and vertical
global coordinates}. The discrete strong form divergence operator, as
given by the incidence matrix $\boldsymbol{\mathsf{E}}^{3,2}$, may itself be decomposed into its horizontal,
$\boldsymbol{\mathsf{E}}^{3,2}_{\parallel}$ and vertical $\boldsymbol{\mathsf{E}}^{3,2}_{\perp}$ components, such that
\begin{equation}
	\nabla\cdot\boldsymbol{u}_{h} \stackrel{\eqref{eq_incidence_matrices}}{=} \gamma_{h}\boldsymbol{\mathsf{E}}^{3,2}\valgebra{\boldsymbol{u}}_{h} =: \gamma_{h}\left(\boldsymbol{\mathsf{E}}^{3,2}_{\parallel} \valgebra{\boldsymbol{v}}_{h} + \boldsymbol{\mathsf{E}}^{3,2}_{\perp} \valgebra{\boldsymbol{w}}_{h}\right)\,.
\end{equation}

The second order HEVI splitting is then given as:
\\
\\
\noindent
\emph{Step 1: Explicit horizontal momentum solve}
\begin{multline}\label{eq::step_1}
\uMassPar{\hat{v}_h}' = \uMassPar{\hat{v}_h}^{n-1} -
2\Delta t\langle\boldsymbol{\beta}_h^{\parallel},q^{\perp,n}_h\times\boldsymbol{\beta}_h^{\parallel}\rangle{\hat{V}_h}^n +
2\Delta t\langle\boldsymbol{\beta}_h^{\parallel}\cdot\boldsymbol{q}_h^{\parallel,n},\beta_h^{\perp}\rangle\hat{W}_h^n + \\
2\Delta t\EdcTPar\qMass{\hat{\Phi}_h}^n +
2\Delta t
\langle\boldsymbol{\beta}_h^{\parallel},\theta_h^n\boldsymbol{\beta}_h^{\parallel}\rangle
\langle\boldsymbol{\beta}_h^{\parallel},\boldsymbol{\beta}_h^{\parallel}\rangle^{-1}
\EdcTPar
\langle\gamma_h,\gamma_h\rangle
{\hat{\Pi}_h}^n
\end{multline}
\noindent
\emph{Step 2: Implicit vertical solve (including horizontal divergence)}
\begin{multline}\label{eq::strang2}
\begin{bmatrix}
\uMassPerp{\hat{w}_h}^{n+1}  \\
\qMass{\hat{\rho}_h}^{n+1}   \\
\qMass{\hat{\Theta}_h}^{n+1} \\
\end{bmatrix} =
\begin{bmatrix}
\uMassPerp{\hat{w}_h}^{n}  \\
\qMass{\hat{\rho}_h}^{n}   \\
\qMass{\hat{\Theta}_h}^{n} \\
\end{bmatrix} - \\
\Delta t
\begin{bmatrix}
%\rMassPerp                                   & -\EdcTPerp\qMass & -\uMassThetaPerp\uMassPerp^{-1}\EdcTPerp\qMass \\
\zero                                        & -\EdcTPerp\qMass & -\uMassThetaPerp\uMassPerp^{-1}\EdcTPerp\qMass \\
\qMass\EdcPerp                               & \zero            & \zero                                          \\
\qMass\EdcPerp\uMassPerp^{-1}\uMassThetaPerp & \zero            & \zero                                          \\
\end{bmatrix}
\begin{bmatrix}
\overline{{\hat{W}}}_h    \\
\overline{{\hat{\Phi}}}_h \\
\overline{{\hat{\Pi}}}_h  \\
\end{bmatrix} - \\
\Delta t
\begin{bmatrix}
\langle\beta_h^{\perp},\boldsymbol{q}_h^{\parallel,n+1/2}\cdot\boldsymbol{\beta}_h^{\parallel}\rangle\overline{\hat{V}}_h \\
\qMass\EdcPar\overline{{\hat{V}}}_h \\
\qMass\EdcPar\uMassPar^{-1}\uMassThetaPar\overline{{\hat{V}}}_h
\end{bmatrix}
\end{multline}
\emph{Step 3: Explicit horizontal momentum solve}
\begin{multline}\label{eq::strang3}
\uMassPar{\hat{v}_h}^{n+1} = \uMassPar{\hat{v}_h}^n -
%\Delta t\rMassPar\overline{{\hat{V}}}_h +
\Delta t\langle\boldsymbol{\beta}_h^{\parallel},q^{\perp,n+1/2}_h\times\boldsymbol{\beta}_h^{\parallel}\rangle\overline{\hat{V}}_h +
\Delta t\langle\boldsymbol{\beta}_h^{\parallel}\cdot\boldsymbol{q}_h^{\parallel,n+1/2},\beta_h^{\perp}\rangle\overline{\hat{W}}_h + \\
\Delta t\EdcTPar\qMass\overline{{\hat{\Phi}}}_h +
\Delta t
\langle\boldsymbol{\beta}_h^{\parallel},\theta_h^{n+1/2}\boldsymbol{\beta}_h^{\parallel}\rangle
\langle\boldsymbol{\beta}_h^{\parallel},\boldsymbol{\beta}_h^{\parallel}\rangle^{-1}
\EdcTPar
\langle\gamma_h,\gamma_h\rangle
\overline{{\hat{\Pi}}}_h
\end{multline}
Left multiplication of \eqref{eq::strang2} by
$[\overline{{\hat{W}}}_h^{\top}\quad\overline{{\hat{\Phi}}}_h^{\top}\quad\overline{{\hat{\Pi}}}_h^{\top}]$
and \eqref{eq::strang3} by $\overline{{\hat{V}}}_h^{\top}$ leads to the cancellation of all forcing terms in both equations,
resulting in the expression
\begin{equation}\label{eq::dHda_dadt}
\overline{{\hat{V}}}_h\uMassPar({\hat{v}}_h^{n+1} - {\hat{v}}_h^{n}) +
\overline{{\hat{W}}}_h\uMassPar({\hat{w}}_h^{n+1} - {\hat{w}}_h^{n}) +
\overline{{\hat{\Phi}}}_h\qMass({\hat{\rho}}_h^{n+1} - {\hat{\rho}}_h^{n}) +
\overline{{\hat{\Pi}}}_h\qMass({\hat{\Theta}}_h^{n+1} - {\hat{\Theta}}_h^{n}) \ne 0
\end{equation}
Note that since ${\boldsymbol{v}_h}'$ and ${\boldsymbol{v}_h}^{n+1}$ differ, energy is not strictly
conserved in the temporal discretisation of the horizontal momentum equation. There is a conservation error
of magnitude $\overline{{\hat{V}}}_h\uMassPar({\hat{v}}_h' - {\hat{v}}_h^{n+1})$, such that 
\eqref{eq::dHda_dadt} is not strictly a discrete analogue of the continuous expression
\begin{equation}
\Bigg(\frac{\delta{H}}{\delta\boldsymbol{y}}\Bigg)^{\top}\cdot
\frac{\partial\boldsymbol{y}}{\partial t} = 0.
\end{equation}
The above formulation is similar to the ``u-forward, pressure-backward'' formulation \cite{Weller13},
by which the horizontal divergence term is treated implicitly (using the most recent horizontal
velocity update), and the horizontal pressure gradient term is treated explicitly. However in the
present case the horizontal pressure gradient term at the final step is evaluated at the same time
level as the (implicit) horizontal divergence term.

%In order to ensure second order accuracy in time, the variational derivatives are given as:
The exact, second order time integrals of the variational derivatives are then given as: 
\begin{subequations}
\begin{align}
\uMassPar\overline{{\hat{V}}}_h =&
\frac{1}{3}\langle\boldsymbol{\beta}_h^{\parallel},\rho_h^n\boldsymbol{\beta}_h^{\parallel}\rangle{\hat{v}_h}^n +
\frac{1}{6}\langle\boldsymbol{\beta}_h^{\parallel},\rho_h^{n+1}\boldsymbol{\beta}_h^{\parallel}\rangle{\hat{v}_h}^n +
\frac{1}{6}\langle\boldsymbol{\beta}_h^{\parallel},\rho_h^n\boldsymbol{\beta}_h^{\parallel}\rangle{\hat{v}_h}' +
\frac{1}{3}\langle\boldsymbol{\beta}_h^{\parallel},\rho_h^{n+1}\boldsymbol{\beta}_h^{\parallel}\rangle{\hat{v}_h}' \label{eq::horiz_mass_flux}\\
\uMassPerp\overline{{\hat{W}}}_h =&
\frac{1}{3}\langle\beta_h^{\perp},\rho_h^n\beta_h^{\perp}\rangle{\hat{w}_h}^n +
\frac{1}{6}\langle\beta_h^{\perp},\rho_h^{n+1}\beta_h^{\perp}\rangle{\hat{w}_h}^n +
\frac{1}{6}\langle\beta_h^{\perp},\rho_h^n\beta_h^{\perp}\rangle{\hat{w}_h}^{n+1} +
\frac{1}{3}\langle\beta_h^{\perp},\rho_h^{n+1}\beta_h^{\perp}\rangle{\hat{w}_h}^{n+1} \\
\qMass\overline{{\hat{\Phi}}}_h =&
\frac{1}{6}\langle\gamma_h,\boldsymbol{v}_h^n\cdot\boldsymbol{\beta}_h^{\parallel}\rangle{\hat{v}_h}^n +
\frac{1}{6}\langle\gamma_h,\boldsymbol{v}'_h\cdot\boldsymbol{\beta}_h^{\parallel}\rangle{\hat{v}_h}^n +
\frac{1}{6}\langle\gamma_h,\boldsymbol{v}'_h\cdot\boldsymbol{\beta}_h^{\parallel}\rangle{\hat{v}_h}' +\notag\\
&
\frac{1}{6}\langle\gamma_h,w_h^n\beta_h^{\perp}\rangle{\hat{w}_h}^n +
\frac{1}{6}\langle\gamma_h,w^{n+1}_h\beta_h^{\perp}\rangle{\hat{w}_h}^n +
\frac{1}{6}\langle\gamma_h,w^{n+1}_h\beta_h^{\perp}\rangle{\hat{w}_h}^{n+1} +
g\langle\gamma_h,\gamma_h\rangle{\hat{z}_h}\\
\qMass\overline{{\hat{\Pi}}}_h =&
\frac{1}{2}\langle\gamma_h,\gamma_h\rangle{\hat{\Pi}_h}^n +
\frac{1}{2}\langle\gamma_h,\gamma_h\rangle{\hat{\Pi}_h}^{n+1}
\end{align}
\end{subequations}
Note that while a second order \blue{``leap-frog'' } method has been used to determine the provisional
velocity, ${\hat{v}_h}'$ in \eqref{eq::step_1}, a first order, forward Euler step is probably also
acceptable, since in either case the resulting mass flux \eqref{eq::horiz_mass_flux} will be second
order in time.

\green{We have not as yet made mention of the spaces used to represent the quantities found within the operator 
in \eqref{eq::strang2}: $\boldsymbol{q}_h$, $\theta_h$. Conservation of energy is satisfied for any representation
these quantities, provided that the skew-symmetry of this operator is preserved. We choose to represent 
$\boldsymbol{q}_h\in\mathcal{W}_h$, such that this may be derived from $\boldsymbol{u}_h$ and $\rho_h$ via a weak 
curl operator \cite{LP20}, while $\theta_h$ is represented in the vertical component of the 
$H(\text{div},\Omega)$ space, $\mathcal{U}_h^{\perp}$ (such that it is $C^0$ continuous in the vertical direction 
only). Previous work has shown that this representation of $\theta_h$ yields an improved representation of the
dispersion relation for linearised buoyancy modes \cite{Melvin18}.}

\green{The implicit vertical system in \eqref{eq::strang2} is solved in serial using a direct LU method, 
while the projections onto the horizontal $\mathcal{U}_h^{\parallel}$ space in \eqref{eq::step_1}, \eqref{eq::strang3}, 
as well as to derive the horizontal pressure gradients and temperature fluxes in \eqref{eq::strang2} and
the horizontal mass flux in \eqref{eq::horiz_mass_flux} are
solved in parallel using the conjugate gradient method. These are both done using the PETSc library }
\cite{petsc-user-ref,petsc-web-page,petsc-efficient}.

\subsection{Stability analysis}

In order to study the stability of the HEVI splitting scheme described above this is applied to the
two dimensional linearised compressible Boussinesq equations \cite{Durran10}, which are given as:
\begin{subequations}
\begin{align}
\frac{\partial u}{\partial t} + \frac{\partial p}{\partial x} &= 0 \\
\frac{\partial w}{\partial t} + \frac{\partial p}{\partial z} - b &= 0 \\
\frac{\partial p}{\partial t} + c^2\Bigg(\frac{\partial u}{\partial x} + \frac{\partial w}{\partial z}\Bigg) &= 0 \\
\frac{\partial b}{\partial t} + N^2w &= 0
\end{align}
\end{subequations}
where $u$ and $w$ are the horizontal and vertical velocities, $p$ is the pressure, $b$ is the buoyancy, $c$ is
the speed of sound and $N$ is the Brunt-V\"ais\"al\"a frequency. Using the splitting scheme detailed above
(but with a first order forward Euler step for the initial provisional velocity, $u'$) and
assuming periodic solutions of the form $(u,w,p,b)(t)e^{ikx + ilz}$ these may be time-stepped as:
\begin{subequations}
\begin{align}
u' &= u^n - \Delta t\ ikp^n \\
w^{n+1} + \frac{\Delta t\ il}{2}p^{n+1} - \frac{\Delta t}{2}b^{n+1} &=
w^{n} - \frac{\Delta t\ il}{2}p^{n} + \frac{\Delta t}{2}b^{n} \\
p^{n+1} + \frac{\Delta t\ ilc^2}{2}w^{n+1} &=
p^{n} - \frac{\Delta t\ ilc^2}{2}w^{n} - \frac{\Delta t\ ikc^2}{2}u^{n} - \frac{\Delta t\ ikc^2}{2}u' \\
b^{n+1} + \frac{\Delta t\ N^2}{2}w^{n+1} &= b^{n} - \frac{\Delta t\ N^2}{2}w^{n} \\
u^{n+1} + \frac{\Delta t\ ik}{2}p^{n+1} &= u^{n} - \frac{\Delta t\ ik}{2}p^{n}
\end{align}
\end{subequations}
\blue{Note that a temporal description of the HEVI scheme that omits the spatial discretisation 
is provided in Appendix A. }
Substituting the expression for $u'$ into that for $p^{n+1}$ and expressing in matrix form gives
\begin{equation}
\begin{bmatrix}
1 & 0                       & \frac{\Delta t\ ik}{2} & 0                   \\
0 & 1                       & \frac{\Delta t\ il}{2} & -\frac{\Delta t}{2} \\
0 & \frac{\Delta t\ ilc^2}{2} & 1                    & 0                   \\
0 & \frac{\Delta t\ N^2}{2}   & 0                    & 1                   \\
\end{bmatrix}
\begin{bmatrix}
u^{n+1} \\ w^{n+1} \\ p^{n+1} \\ b^{n+1} \\
\end{bmatrix} =
\begin{bmatrix}
1              & 0                        & -\frac{\Delta t\ ik}{2}          & 0                   \\
0              & 1                        & -\frac{\Delta t\ il}{2}          & \frac{\Delta t}{2}  \\
-\Delta t\ ikc^2 & -\frac{\Delta t\ ilc^2}{2} & 1 - \frac{\Delta t^2k^2c^2}{2} & 0                   \\
0              & -\frac{\Delta t\ N^2}{2}   & 0                              & 1                   \\
\end{bmatrix}
\begin{bmatrix}
u^{n} \\ w^{n} \\ p^{n} \\ b^{n} \\
\end{bmatrix}.
\end{equation}

In order to analyse the stability of the splitting scheme we inspect the eigenvalues of the operator
\begin{equation}
\boldsymbol{\mathsf{A}} =
\begin{bmatrix}
1 & 0                       & \frac{\Delta t\ ik}{2} & 0                   \\
0 & 1                       & \frac{\Delta t\ il}{2} & -\frac{\Delta t}{2} \\
0 & \frac{\Delta t\ ilc^2}{2} & 1                    & 0                   \\
0 & \frac{\Delta t\ N^2}{2}   & 0                    & 1                   \\
\end{bmatrix}^{-1}
\begin{bmatrix}
1              & 0                        & -\frac{\Delta t\ ik}{2}          & 0                   \\
0              & 1                        & -\frac{\Delta t\ il}{2}          & \frac{\Delta t}{2}  \\
-\Delta t\ ikc^2 & -\frac{\Delta t\ ilc^2}{2} & 1 - \frac{\Delta t^2k^2c^2}{2} & 0                   \\
0              & -\frac{\Delta t\ N^2}{2}   & 0                              & 1                   \\
\end{bmatrix}
\end{equation}
For a given eigenvalue, $\lambda_i$, this will be unstable \blue{if } $\lambda_i\lambda_i^* > 1$. The four eigenvalues
correspond to two acoustic modes, which differ in the sign of their imaginary component only, and two
gravity (buoyancy) modes which similarly differ only in the sign of the imaginary term.

The results for the new HEVI scheme are compared to the standard second order trapezoidal HEVI scheme,
for which the vertical terms are all treated implicitly, and the horizontal terms explicitly. This scheme
may be written as:
\begin{equation}
\begin{bmatrix}
1 & 0                       & 0                    & 0                   \\
0 & 1                       & \frac{\Delta t\ il}{2} & -\frac{\Delta t}{2} \\
0 & \frac{\Delta t\ ilc^2}{2} & 1                    & 0                   \\
0 & \frac{\Delta t\ N^2}{2}   & 0                    & 1                   \\
\end{bmatrix}
\begin{bmatrix}
u' \\ w' \\ p' \\ b' \\
\end{bmatrix} =
\begin{bmatrix}
1              & 0                        & -\Delta t\ ik           & 0                  \\
0              & 1                        & -\frac{\Delta t\ il}{2} & \frac{\Delta t}{2} \\
-\Delta t\ ikc^2 & -\frac{\Delta t\ ilc^2}{2} & 1                     & 0                  \\
0              & \frac{\Delta t\ N^2}{2}    & 0                     & 1                  \\
\end{bmatrix}
\begin{bmatrix}
u^n \\ w^n \\ p^n \\ b^n \\
\end{bmatrix}
\end{equation}

\begin{equation}
\begin{bmatrix}
u^{n+1} \\ w^{n+1} \\ p^{n+1} \\ b^{n+1} \\
\end{bmatrix} =
\begin{bmatrix}
1                        & 0                        & -\frac{\Delta t\ ik}{2} & 0                  \\
0                        & 1                        & -\frac{\Delta t\ il}{2} & \frac{\Delta t}{2} \\
-\frac{\Delta t\ ikc^2}{2} & -\frac{\Delta t\ ilc^2}{2} & 1                     & 0                  \\
0                        & -\frac{\Delta t\ N^2}{2}   & 0                     & 1                  \\
\end{bmatrix}
\begin{bmatrix}
u^{n} \\ w^{n} \\ p^{n} \\ b^{n} \\
\end{bmatrix} +
\begin{bmatrix}
0                        & 0                        & -\frac{\Delta t\ ik}{2} & 0                  \\
0                        & 0                        & -\frac{\Delta t\ il}{2} & \frac{\Delta t}{2} \\
-\frac{\Delta t\ ikc^2}{2} & -\frac{\Delta t\ ilc^2}{2} & 0                     & 0                  \\
0                        & -\frac{\Delta t\ N^2}{2}   & 0                     & 0                  \\
\end{bmatrix}
\begin{bmatrix}
u' \\ w' \\ p' \\ b' \\
\end{bmatrix}
\end{equation}
Combining the above expressions gives a linear operator by which the solution at time level $n+1$
may be determined from the solution at time level $n$, from which the corresponding eigenvalues may
be derived.

For completeness we also compare against the Crank-Nicolson scheme, which is semi-implicit
in both horizontal and vertical dimensions, and neutrally stable for all vertical and horizontal
wave numbers for both buoyancy and acoustic modes. The amplification factors for the
Crank-Nicolson scheme are given as:

\begin{equation}
\boldsymbol{\mathsf{A}} =
\begin{bmatrix}
1 & 0                       & \frac{\Delta t\ ik}{2} & 0                   \\
0 & 1                       & \frac{\Delta t\ il}{2} & -\frac{\Delta t}{2} \\
\frac{\Delta t\ ikc^2}{2} & \frac{\Delta t\ ilc^2}{2} & 1                    & 0                   \\
0 & \frac{\Delta t\ N^2}{2}   & 0                    & 1                   \\
\end{bmatrix}^{-1}
\begin{bmatrix}
1              & 0                        & -\frac{\Delta t\ ik}{2}          & 0                   \\
0              & 1                        & -\frac{\Delta t\ il}{2}          & \frac{\Delta t}{2}  \\
-\frac{\Delta t\ ikc^2}{2} & -\frac{\Delta t\ ilc^2}{2} & 1 & 0                   \\
0              & -\frac{\Delta t\ N^2}{2}   & 0                              & 1                   \\
\end{bmatrix}
\end{equation}

\begin{figure}[!hbtp]
\begin{center}
\includegraphics[width=0.32\textwidth,height=0.24\textwidth]{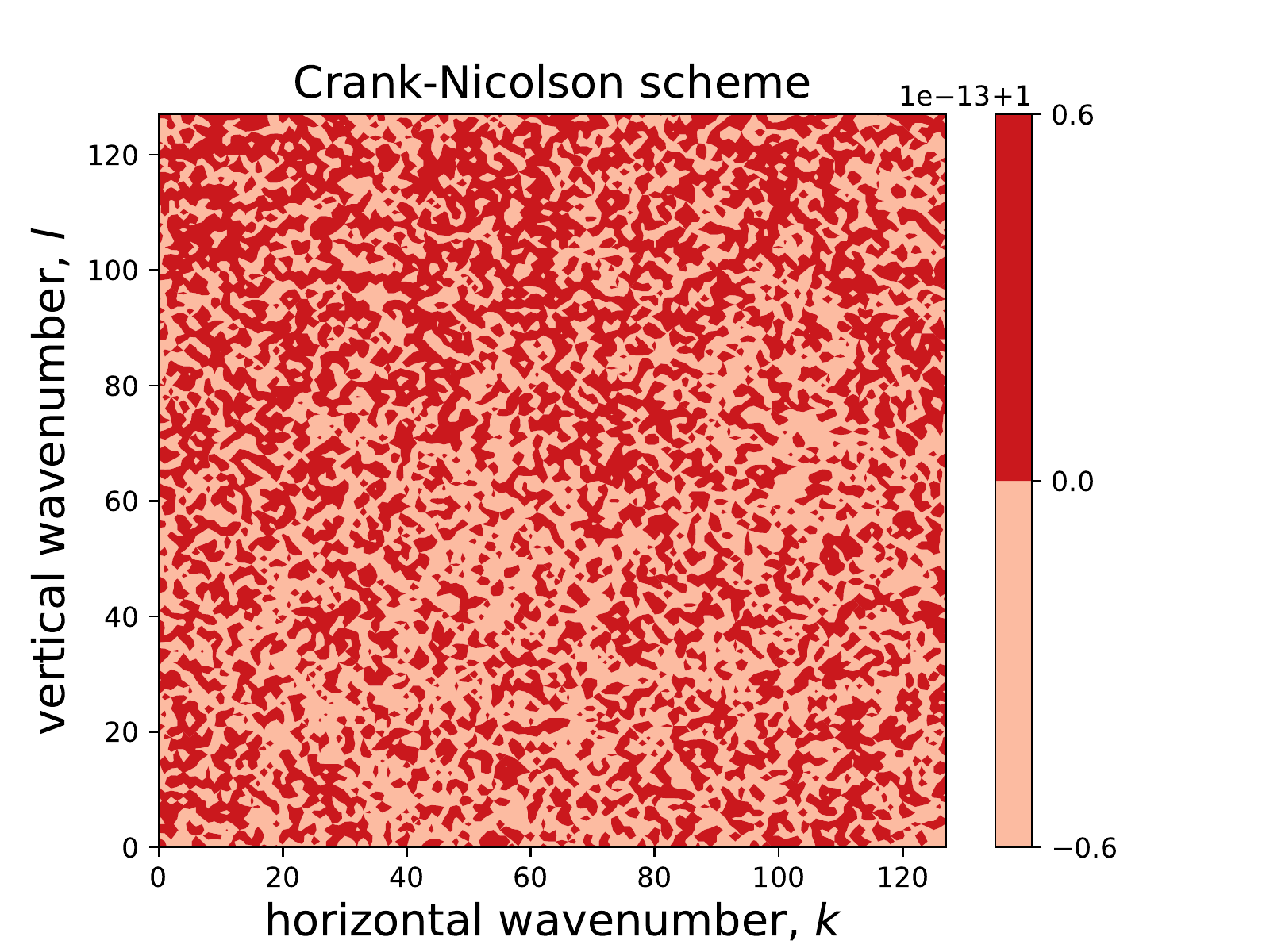}
\includegraphics[width=0.32\textwidth,height=0.24\textwidth]{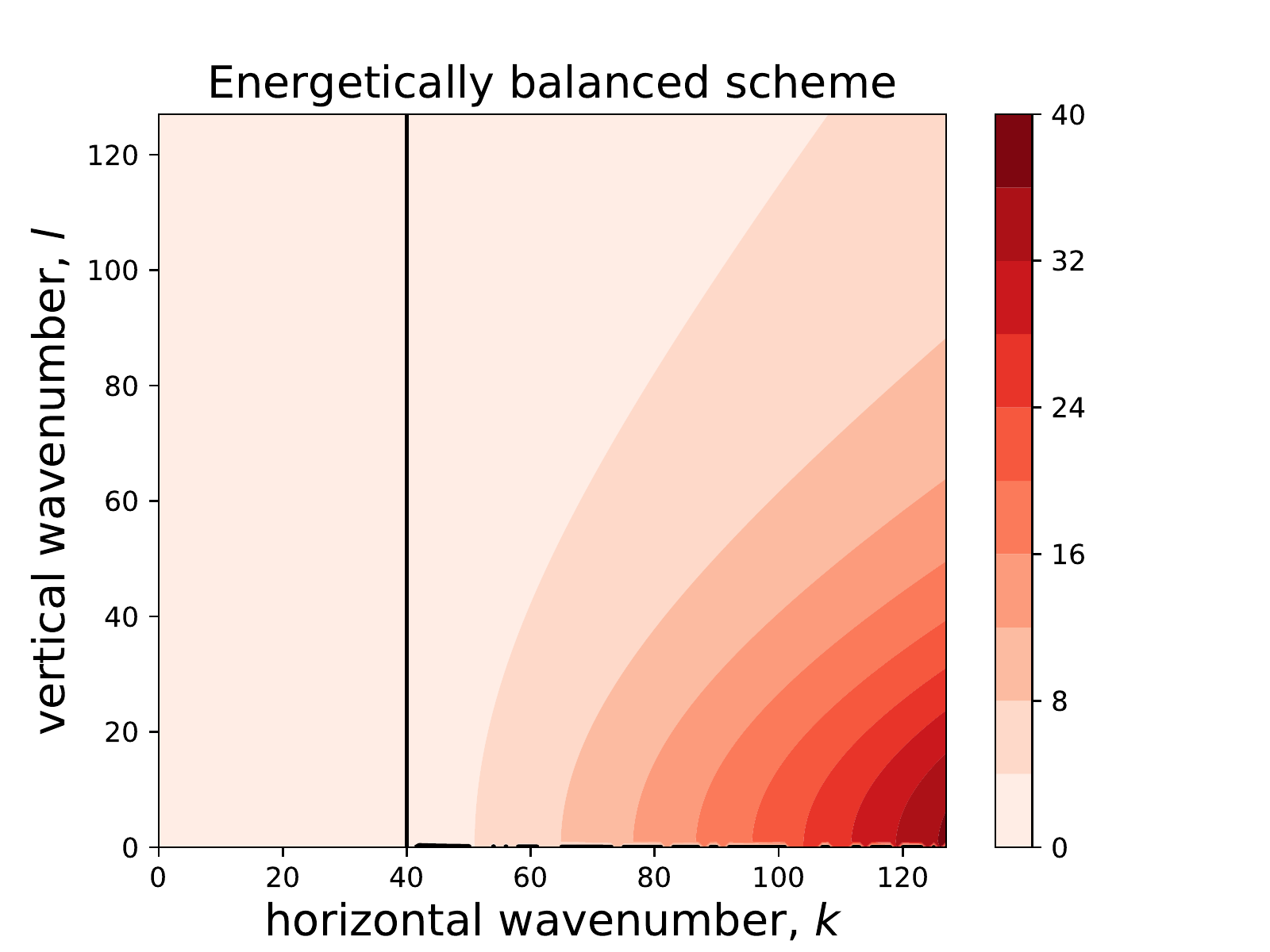}
\includegraphics[width=0.32\textwidth,height=0.24\textwidth]{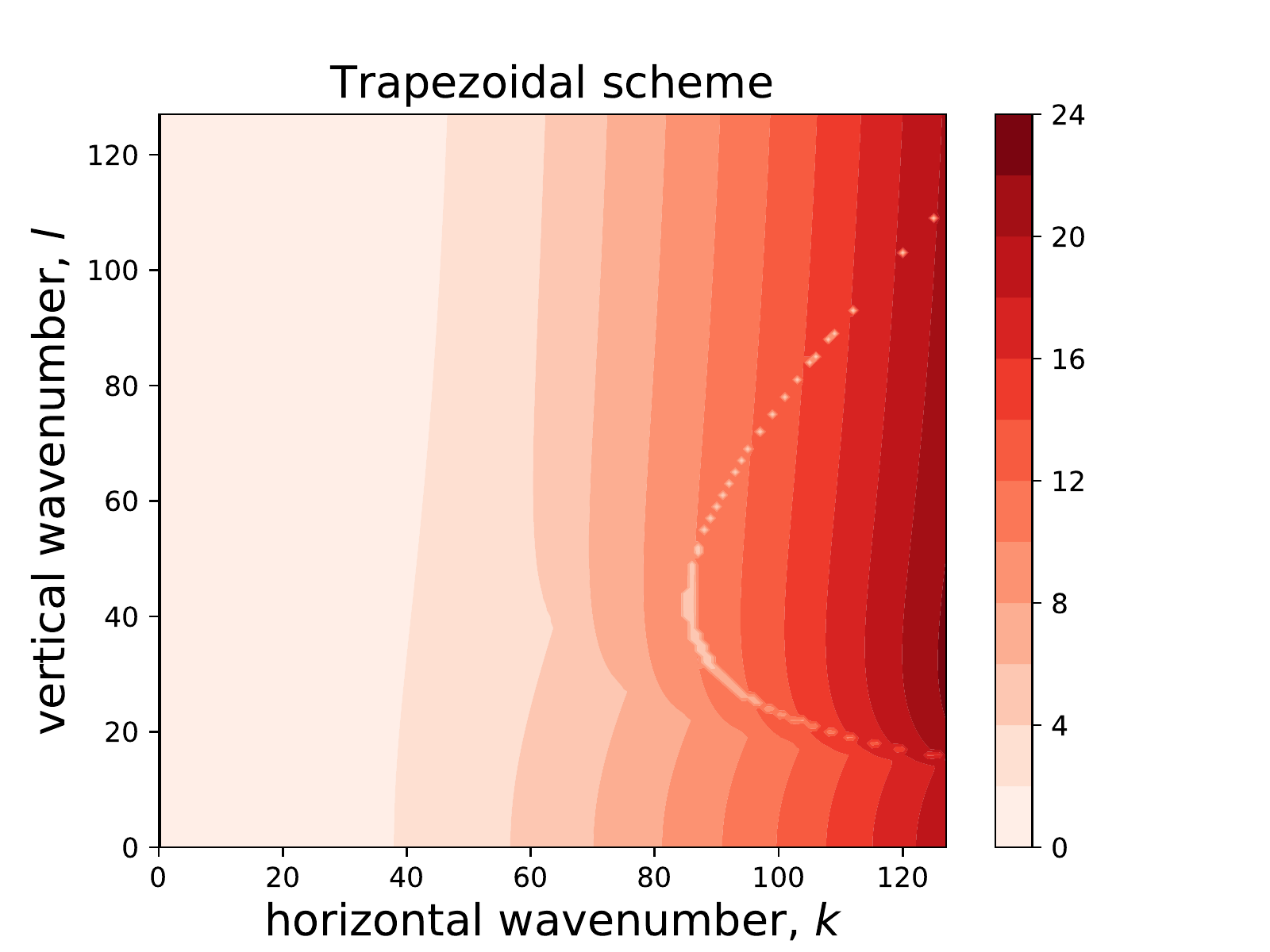}
\caption{Amplification factor for the acoustic waves using the Crank-Nicolson (left), new (center) and
trapezoidal HEVI (right) splitting. The black line on the center plot indicates an amplification factor
of $1.0 + 10^{-12}$, indicating the onset of unstable modes. \green{The Crank-Nicolson scheme is neutrally
stable for all modes, the new scheme is neutrally stable for all modes below a finite horizontal wave
number and the trapezoidal scheme is unstable for all modes.}}
\label{fig::accoustic}
\end{center}
\end{figure}

\begin{figure}[!hbtp]
\begin{center}
\includegraphics[width=0.32\textwidth,height=0.24\textwidth]{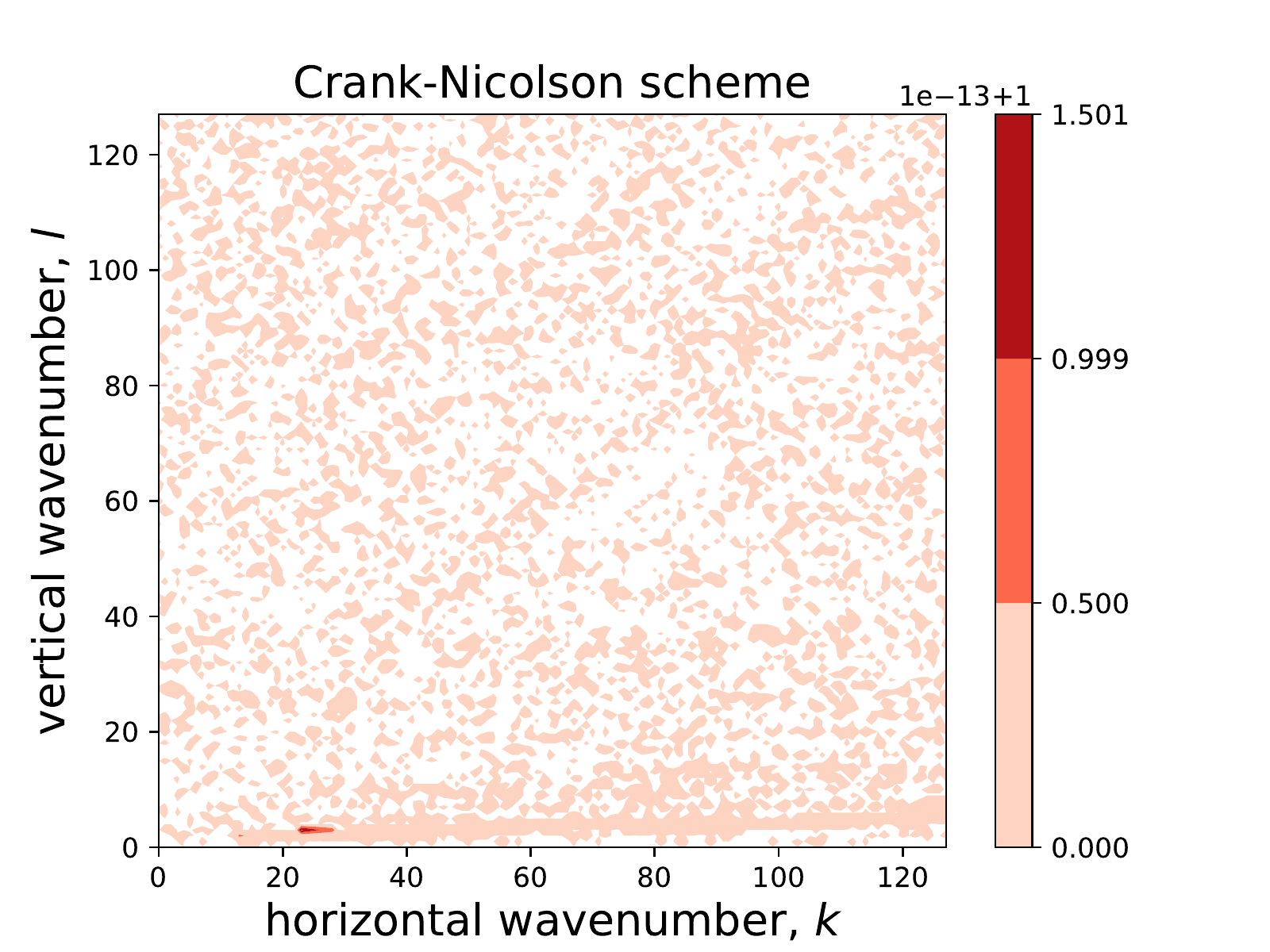}
\includegraphics[width=0.32\textwidth,height=0.24\textwidth]{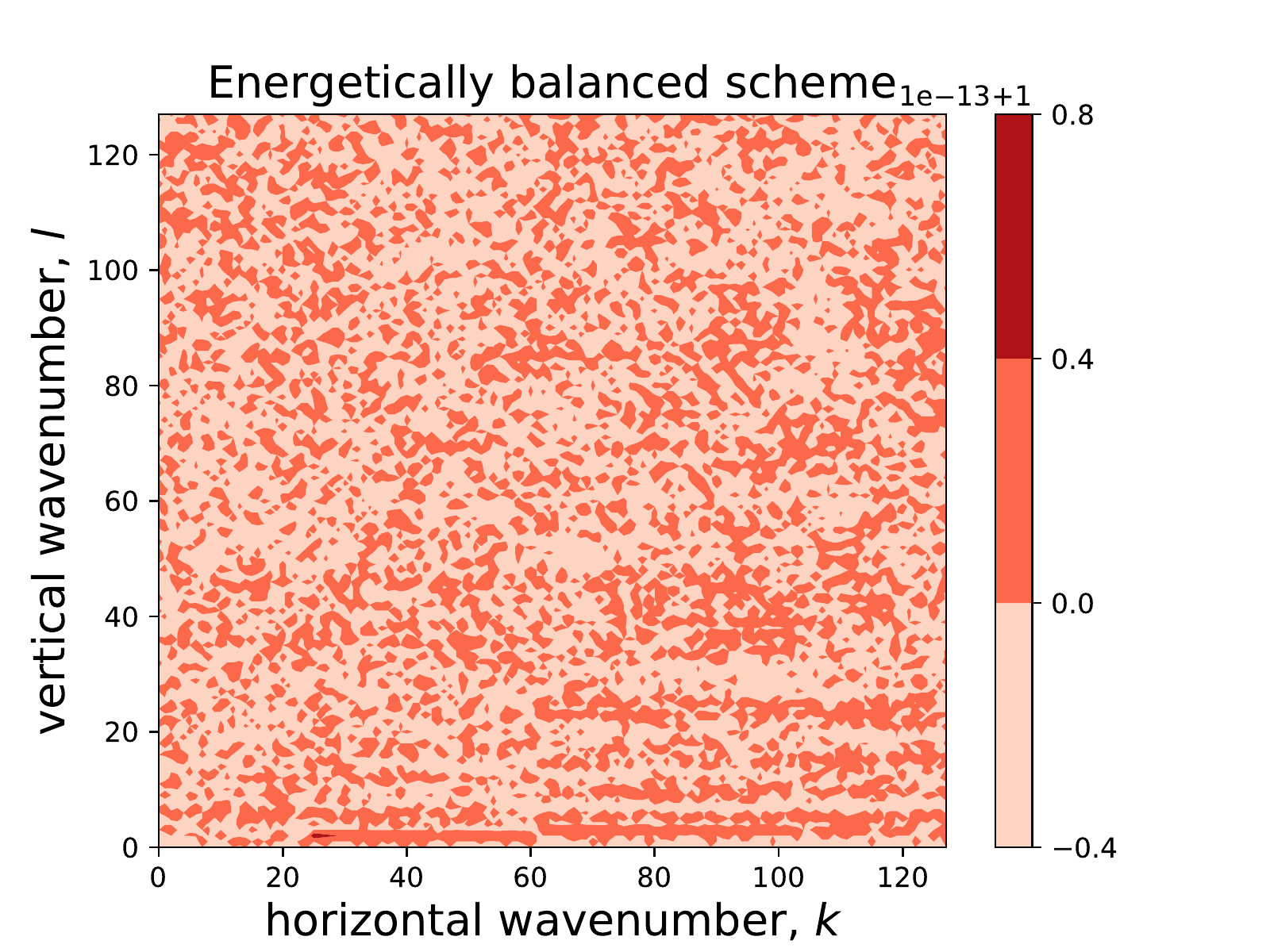}
\includegraphics[width=0.32\textwidth,height=0.24\textwidth]{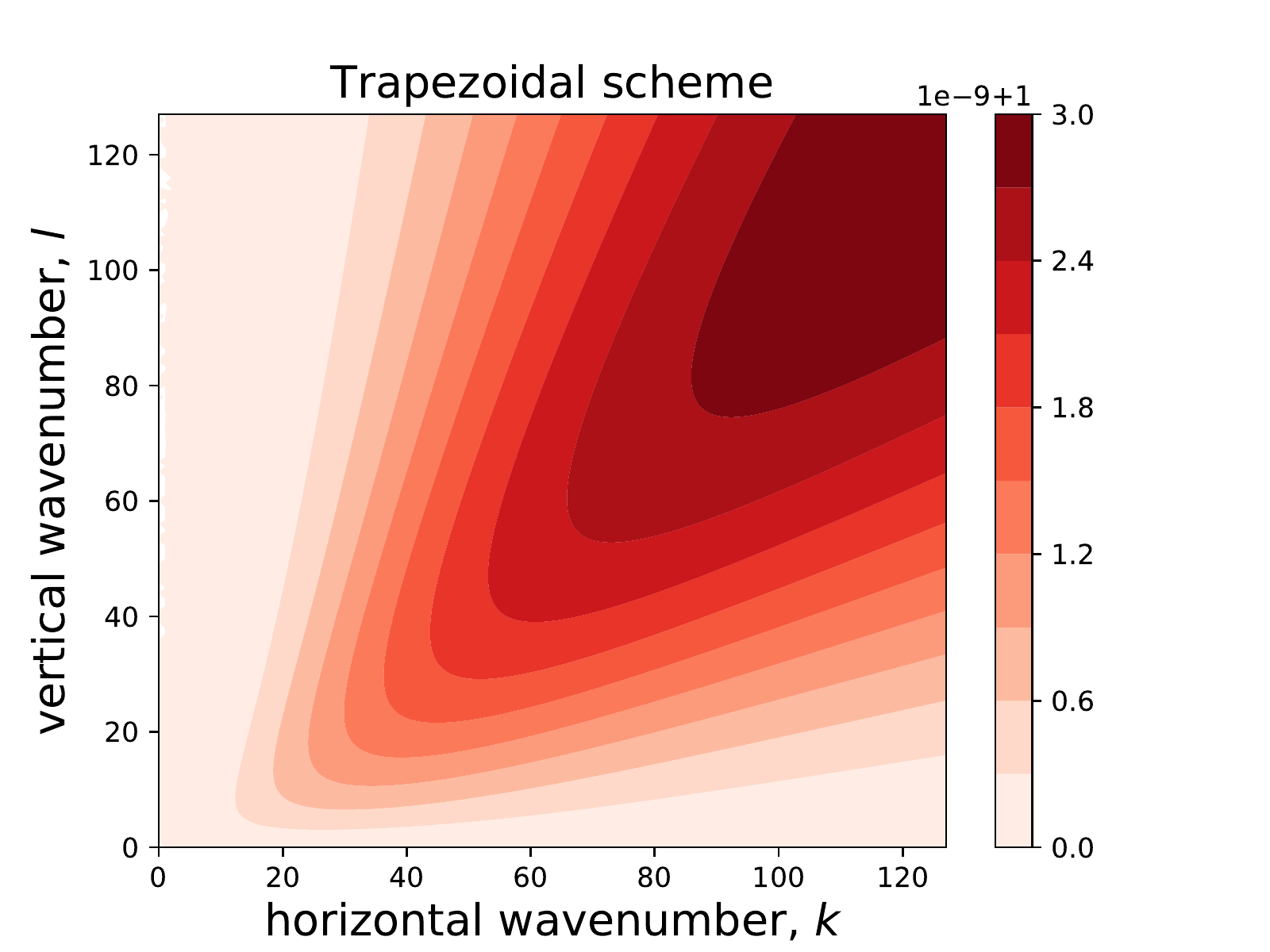}
\caption{Amplification factor for the gravity waves using the Crank-Nicolson (left), new (center) and
trapezoidal HEVI (right) splitting.
\green{The Crank-Nicolson scheme and the new scheme are both neutrally
stable for all modes, while the trapezoidal scheme is unstable for all modes.}}
\label{fig::gravity}
\end{center}
\end{figure}

The amplification factors for the acoustic and gravity modes for the three different time splitting
schemes are given in figs. \ref{fig::accoustic} and \ref{fig::gravity} respectively for a periodic domain of
size 1000.0m in both dimensions, with $c=340.0\mathrm{ms}^{-1}$, $N=0.01\mathrm{s}^{-1}$. While the
Crank-Nicolson scheme is neutrally stable, $|\lambda| = 1$ for both modes, and the trapezoidal HEVI
scheme is unstable, $|\lambda| > 1$ for both modes
(albeit with a very small amplification factor for the gravity modes), the energetically balanced 
splitting scheme is
neutrally stable for all acoustic modes below some horizontal wave number, and unconditionally
stable for all gravity modes. This represents a qualitative improvement on the standard second order
HEVI splitting.

\section{Potential temperature upwinding}\label{sec::upwind}

This section describes the suppression of spurious oscillations in potential temperature via
upwinding of the potential temperature diagnostic equation in an energetically consistent manner.
Since the potential vorticity and potential temperature are represented within the skew-symmetric
operator for the discrete variational form of the compressible Euler equations \eqref{eq::de},
these quantities may be upwinded so as to suppress oscillations without breaking energy conservation
(so long as skew-symmetry is preserved). This strategy has been used previously to suppress
oscillations for the shallow water equations by upwinding potential vorticity so as to dissipate
potential enstrophy \cite{SB85,MC14,Lee20b}, and for the compressible Euler equations
by upwinding potential temperature so as to generate entropy \cite{WCB20}. In the present context
the upwinding of potential temperature is applied not to the material form advection equation,
but rather to the diagnostic equation by which the potential temperature is derived from the
prognostic variables of density \eqref{eq::ce_rho} and density-weighted potential temperature (1c).

While this is typically achieved for variational methods by augmenting the test space with a linear
correction which allows for dissipation in the direction of the flow \cite{BH82}, in the present
case this is achieved by integrating the test functions backwards along velocity characteristics \cite{Lee20b}.
Here the upwinding is applied in the vertical direction only. Since the discrete potential
temperature is represented as $\theta_h\in\mathcal{U}_h^{\perp}$ \cite{LP20}, for which the vertical
component of the tensor product basis is represented by piecewise linear polynomials, for this case
the upwinded bases are effectively equivalent to the SUPG method \cite{BH82,WCB20}.
\blue{However in the present case }
this upwinding is applied to the diagnostic equation \blue{(the discrete form of 
the relation $\theta:=\Theta/\rho$, \eqref{eq::diag_theta})}, and not the material form
of the advection equation, \blue{$\partial\theta/\partial t + \boldsymbol{u}\cdot\nabla\theta = 0$},
\blue{as is customary in the SUPG method. Since $C^0$ continuity for the $\mathcal{U}_h^{\perp}$
space is enforced only in the vertical direction, upwinding of $\theta_h$ in the horizontal directions 
will have no effect.}

The upwinded trial functions, $\beta_h^{\perp,u}\in\mathcal{U}_h^{\perp}$ are computed by evaluating
the quadrature points at vertically \emph{downwind} locations, $\zeta^d$, which are computed to first
order as \cite{Lee20b}:
\begin{equation}
\zeta^d = \zeta + 0.5\Delta t\sum_i\hat{w}_i\beta_i^{\perp}(\xi,\eta,\zeta),
\end{equation}
where $\xi$, $\eta$ and $\zeta$ are the local element coordinates. Note that in the above expression
the velocity is interpolated to its local element value, and not its global physical value.
Note also that the choice to integrate the vertical quadrature points downwind for a time of $0.5\Delta t$
is arbitrary, and may be tuned to optimise the balance between reduced oscillation and less diffusive
solutions.
The upwinded trial functions are then evaluated as $\beta_h^{\perp,u} = \beta_h^{\perp}(\xi,\eta,\zeta^d)$.
The upwinded potential temperature is then diagnosed as:
\begin{equation}\label{eq::diag_theta}
\langle\beta_h^{\perp,u},\rho_h\beta_h^{\perp}\rangle\hat{\theta}_h = \langle\beta_h^{\perp,u},\gamma_h\rangle\hat{\Theta}_h.
\end{equation}
Note that in the above expression only the test functions for the potential temperature diagnostic
equation are evaluated at upstream locations, while the trial functions remain static. Moreover
since continuity is enforced only in the vertical direction for the trial space of the potential
temperature, $\mathcal{U}_h^{\perp}$, it is only in this direction that upwinding has any material
effect.

While the variational form of the compressible Euler equations \eqref{eq::ce_var} conserve entropy,
$s = c_p\log(\theta) + s_0$
(see Appendix B), this is not true in the discrete form, since transcendental functions such as logarithms
cannot be integrated exactly in the discrete form. As such, while the upwinding
of the test functions as described above will generate additional entropy, there are already additional
entropy sources at the discrete level.

\section{Results}\label{sec::test}

\subsection{Baroclinic instability}

\begin{figure}[!hbtp]
\begin{center}
\includegraphics[width=0.48\textwidth,height=0.36\textwidth]{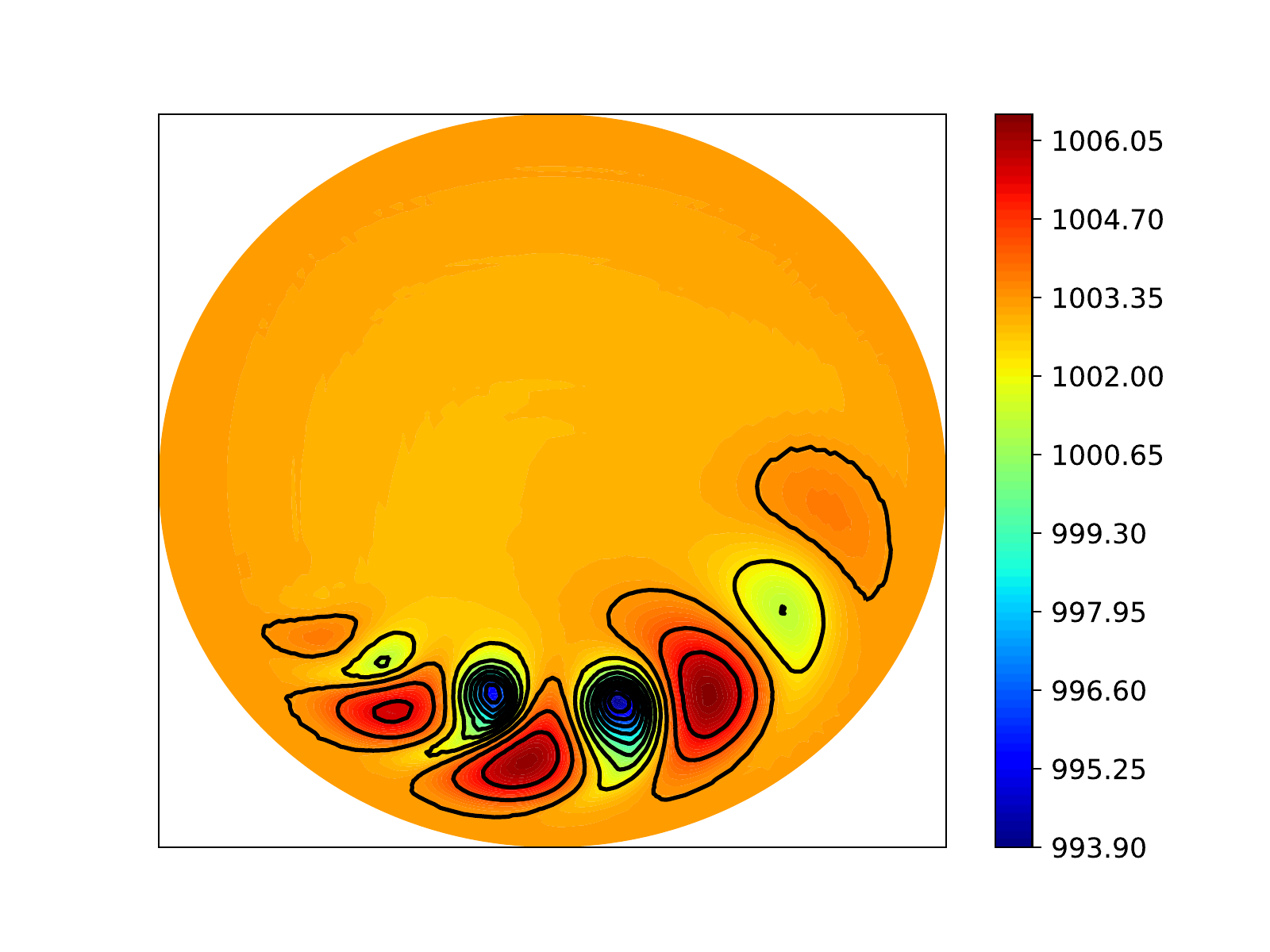}
\includegraphics[width=0.48\textwidth,height=0.36\textwidth]{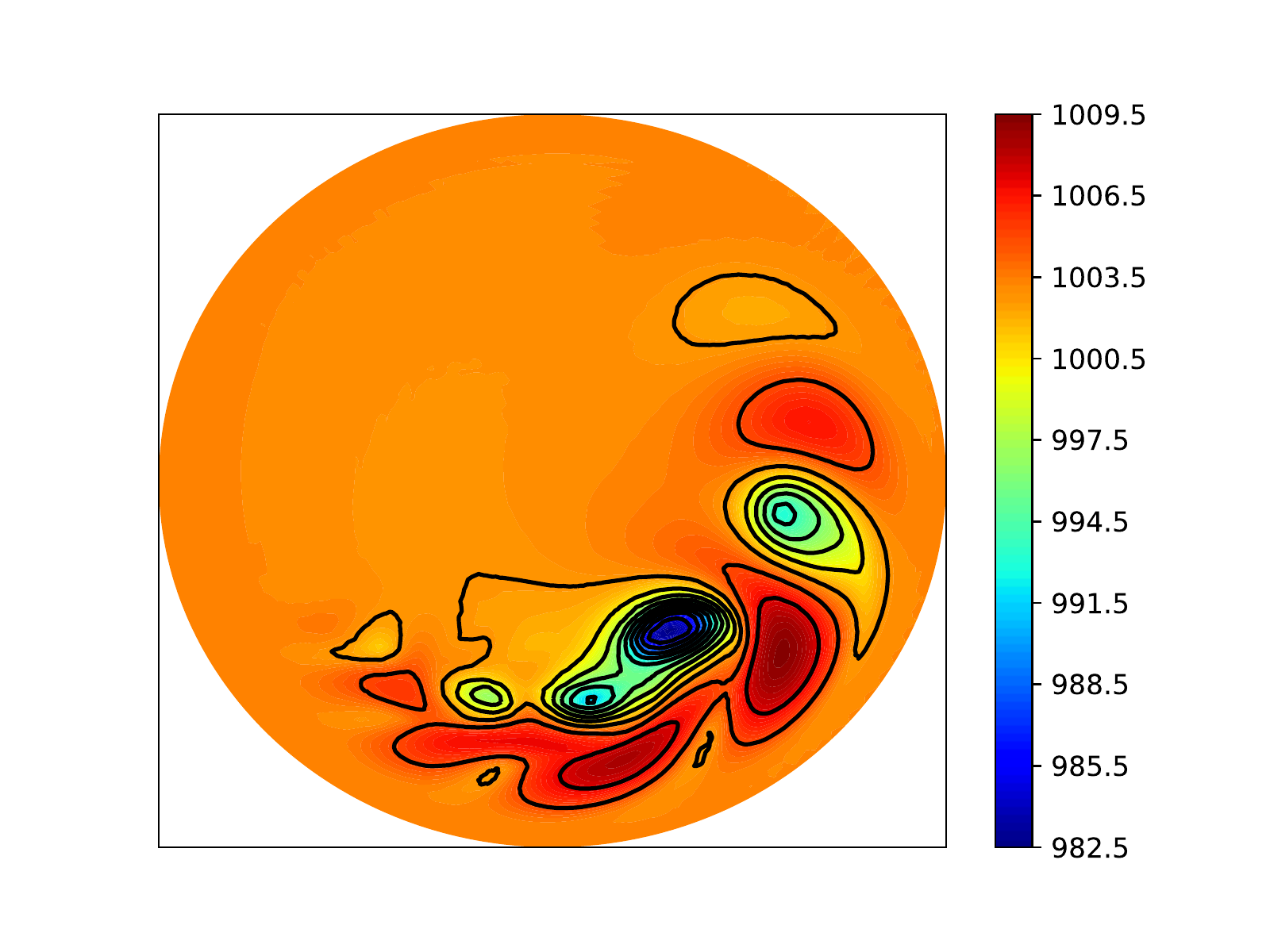}
\caption{\green{Baroclinic instability test case: } Bottom level Exner pressure, days 8 and 10.}
\label{fig::exner_1}
\end{center}
\end{figure}

\begin{figure}[!hbtp]
\begin{center}
\includegraphics[width=0.48\textwidth,height=0.36\textwidth]{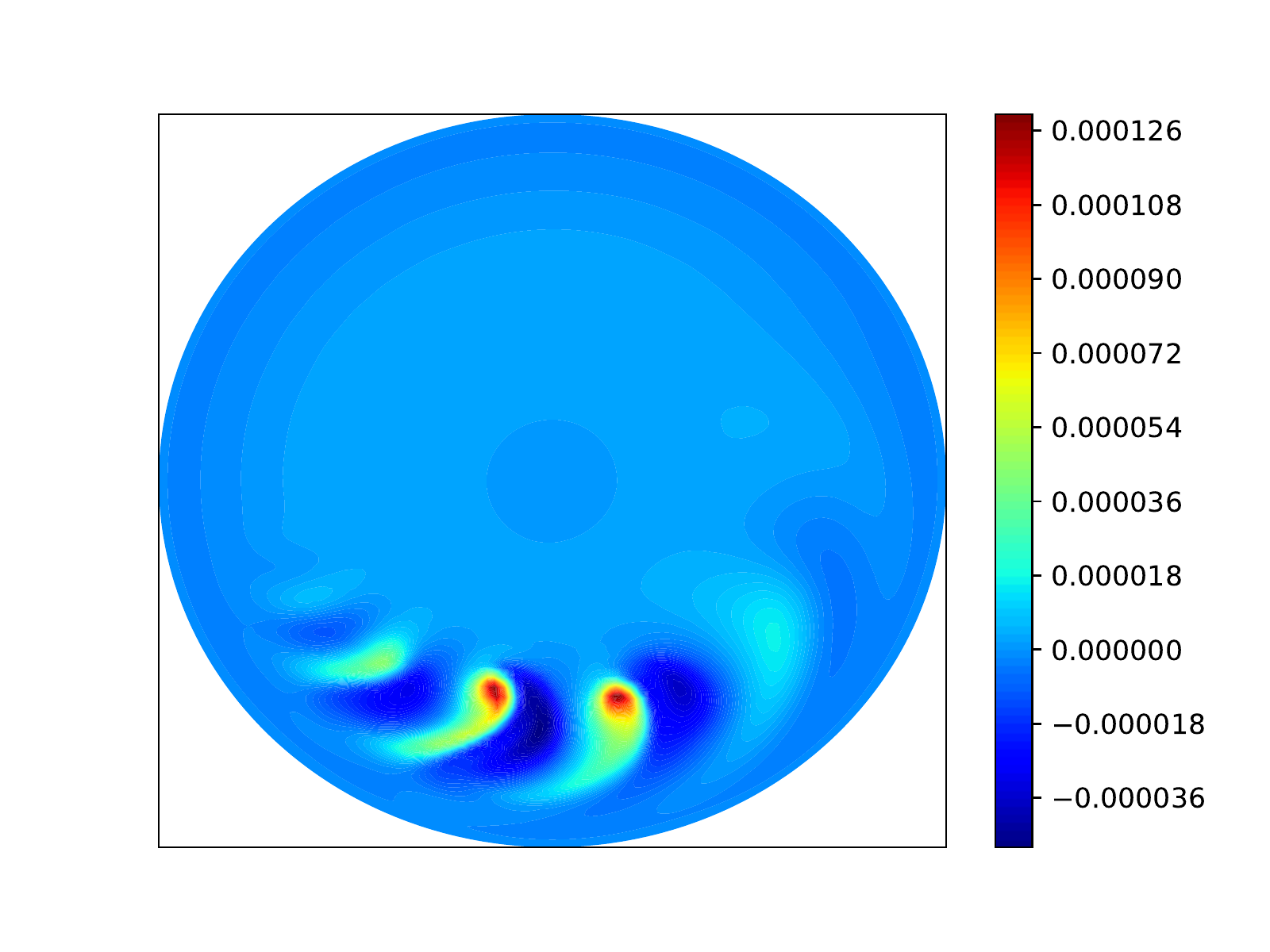}
\includegraphics[width=0.48\textwidth,height=0.36\textwidth]{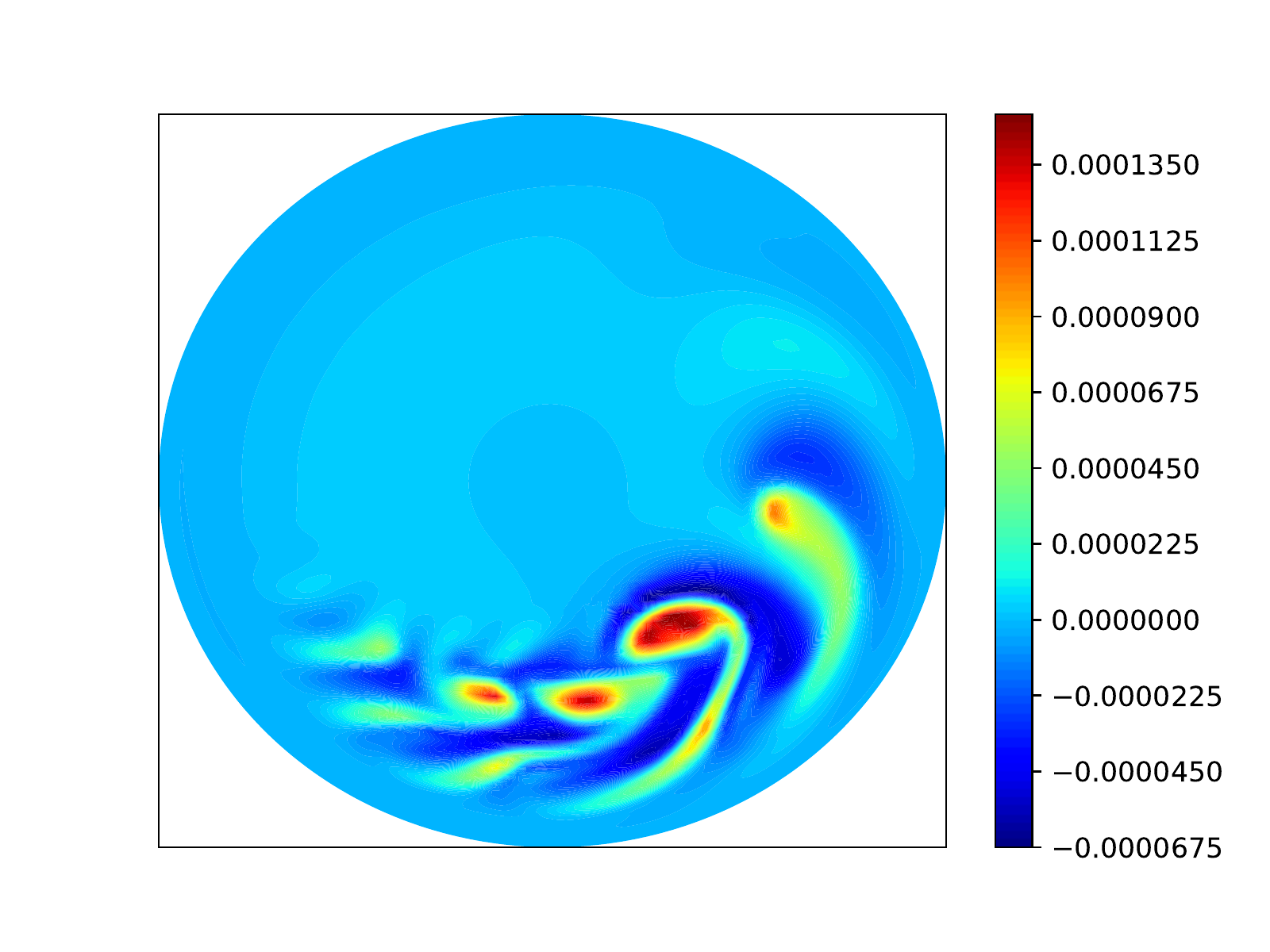}
\caption{\green{Baroclinic instability test case: } Vertical vorticity component at $z=1.57km$, days 8 and 10.}
\label{fig::vorticity_1}
\end{center}
\end{figure}

\begin{figure}[!hbtp]
\begin{center}
\includegraphics[width=0.48\textwidth,height=0.36\textwidth]{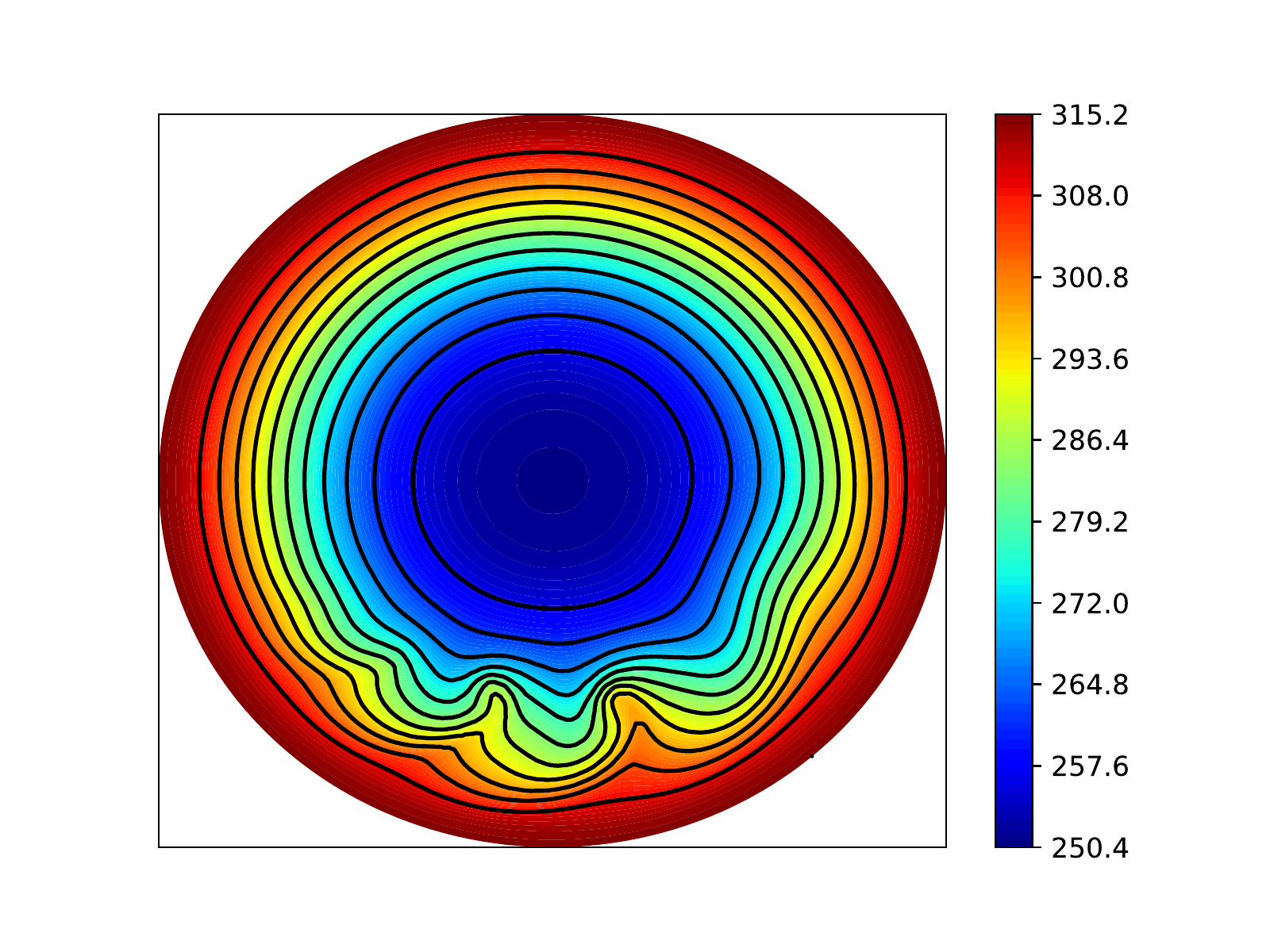}
\includegraphics[width=0.48\textwidth,height=0.36\textwidth]{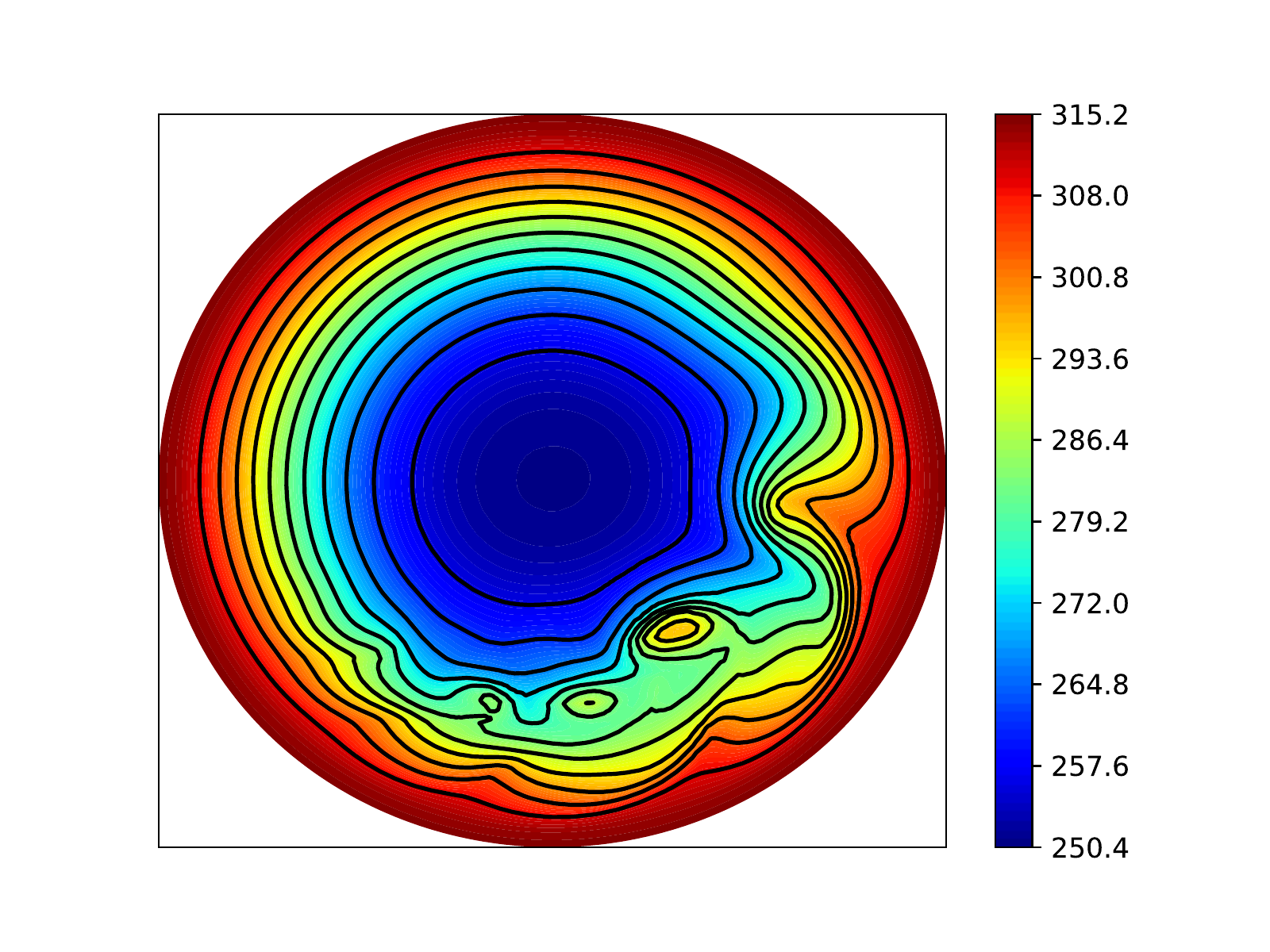}
\caption{\green{Baroclinic instability test case: } Potential temperature at $z=1.57km$, days 8 and 10.}
\label{fig::theta_1}
\end{center}
\end{figure}

\begin{figure}[!hbtp]
\begin{center}
\includegraphics[width=0.48\textwidth,height=0.36\textwidth]{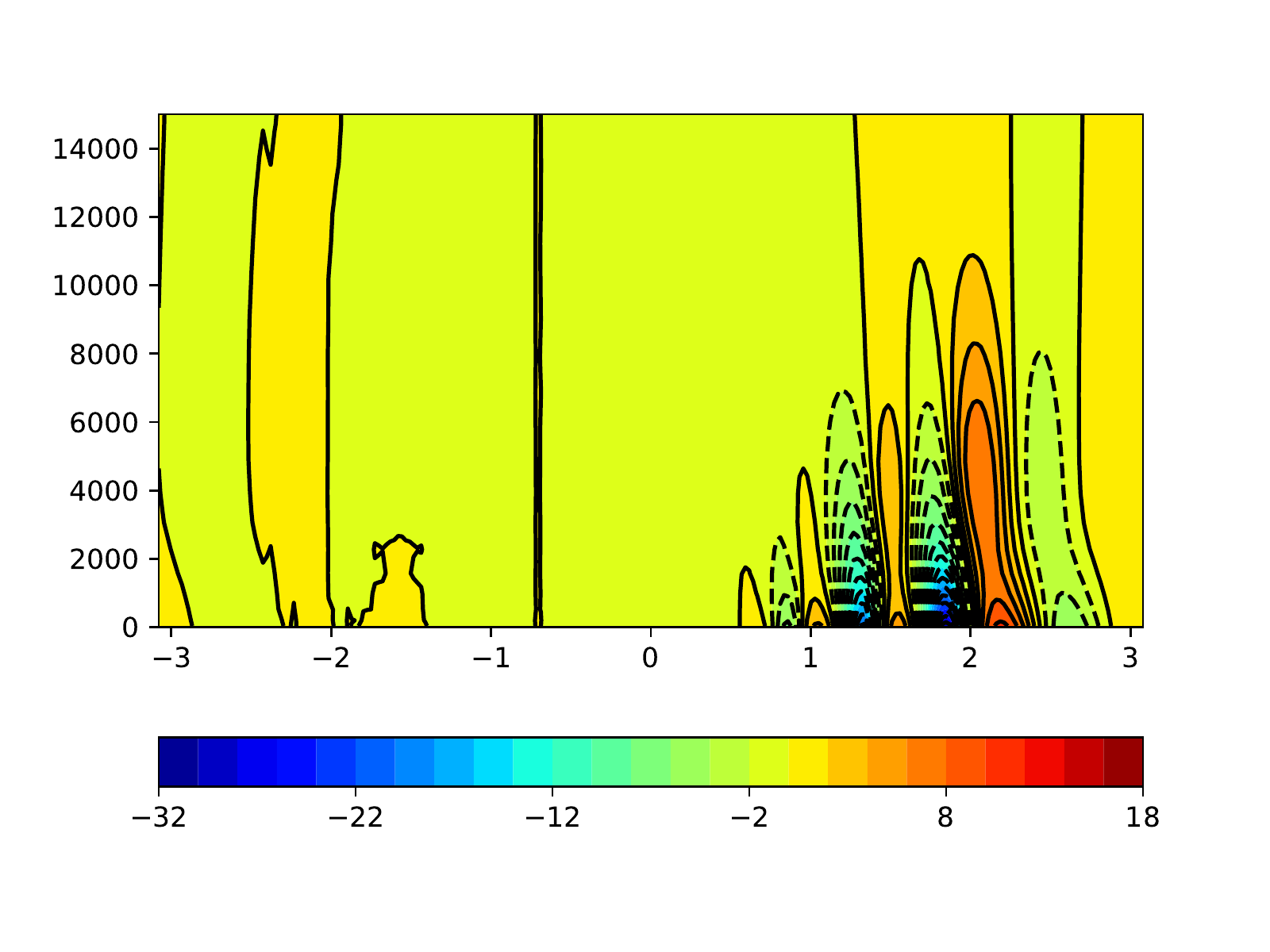}
\includegraphics[width=0.48\textwidth,height=0.36\textwidth]{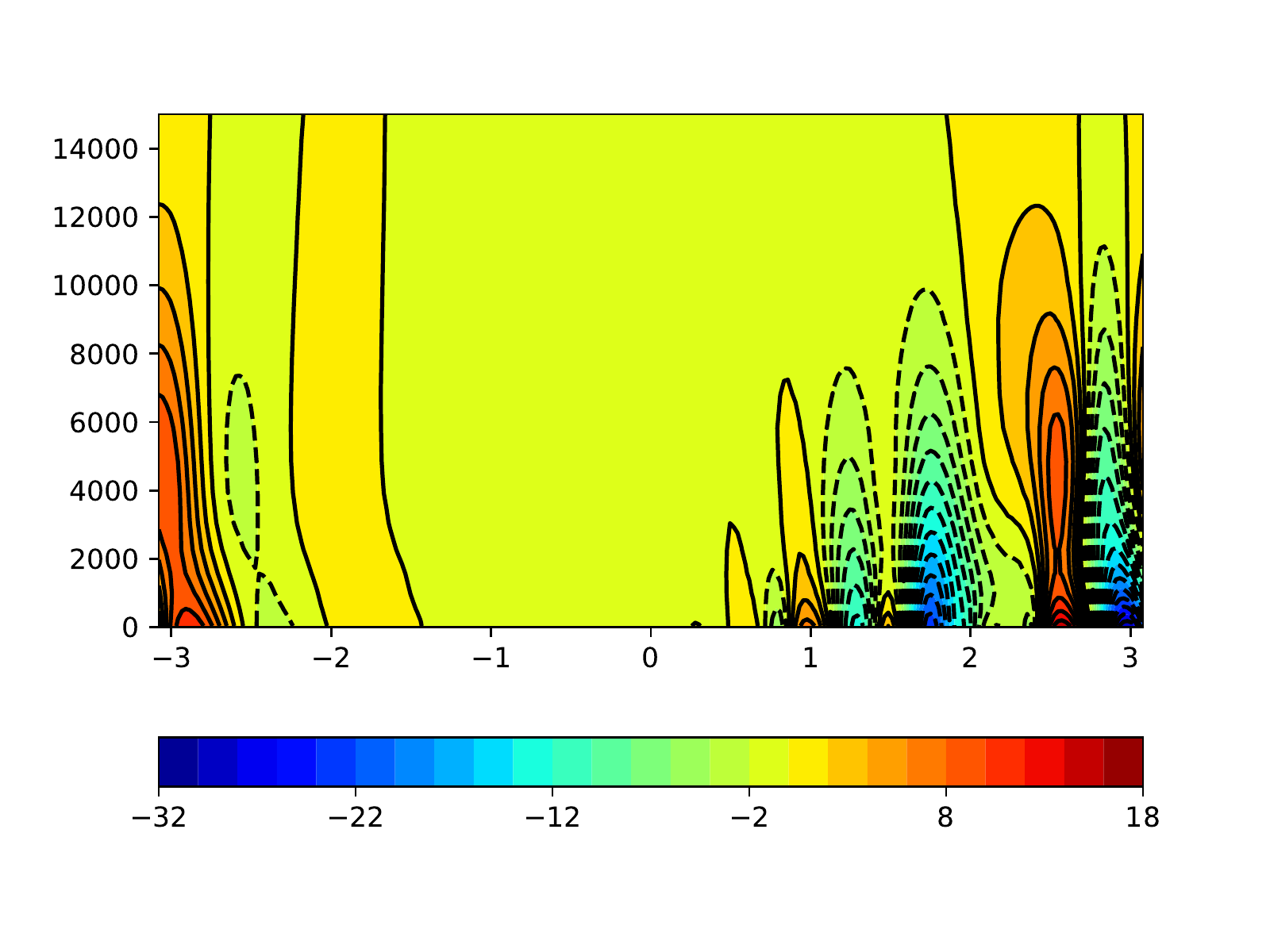}
\caption{\green{Baroclinic instability test case: } Zonal-vertical pressure slice at $50^{\circ}N$, days 8 and 10.}
\label{fig::pressure}
\end{center}
\end{figure}

In order to verify the energetically balanced HEVI scheme, this is applied to the solution of a
standard test case for baroclinic instability on the sphere \cite{UMJS14}. This test involves a
small perturbation to an otherwise hydrostatically and geostrophically balanced atmosphere, which
over several days generates a baroclinic wave. \blue{The model is stabilised using a biharmonic
viscosity term on both the horizontal potential temperature equation and momentum equation (which 
does not conserve energy due to the explicit time stepping), and a Rayleigh friction term 
in the top layers of the atmosphere in the vertical (which is not strictly necessary however helps 
to accelerate the initial hydrostatic adjustment process) \cite{LP20}. }
Figures \ref{fig::exner_1}, \ref{fig::vorticity_1}
and \ref{fig::theta_1} show the surface level Exner pressure, as well as the vertical vorticity
component and potential temperature at a height of $z=1.57\mathrm{km}$ at days 8 and 10
respectively. \blue{As an initial validation exercise } these compare well against previously 
published results \cite{UMJS14,Lee20}, 
\blue{however they do not allow us to draw any particular conclusion as to the superiority
of any particular scheme}.
The signature of the baroclinic instability is also observed in the perturbation to the
steady hydrostatic pressure profile, for which a vertical cross section at $50^{\circ}\mathrm{N}$
is shown at days 8 and 10 in fig \ref{fig::pressure}.

\subsubsection{Comparison to TRAP(2,3,2)}

The energetically balanced HEVI scheme is compared to the horizontally third order, vertically
second order trapezoidal scheme, TRAP(2,3,2) \cite{Weller13,Lock14}. In most respects the energetics
of the two schemes are similar, despite the fact that the TRAP(2,3,2) scheme is not energetically
balanced, but is higher order in the horizontal (also a time step of $\Delta t = 60\mathrm{s}$ is
used for the energetically balanced scheme, while the TRAP(2,3,2) simulation was run with a time
step of $\Delta t = 120\mathrm{s}$). 
\blue{While the TRAP(2,3,2) scheme is stable for longer time steps than the energetically balanced 
scheme, it also involves an additional substep in both the vertical and horizontal directions. 
Both schemes are stabilised using biharmonic viscosity 
in the horizontal and Rayleigh friction in the top layers in the vertical \cite{LP20,Lee20}. }

\begin{figure}[!hbtp]
\begin{center}
\includegraphics[width=0.48\textwidth,height=0.36\textwidth]{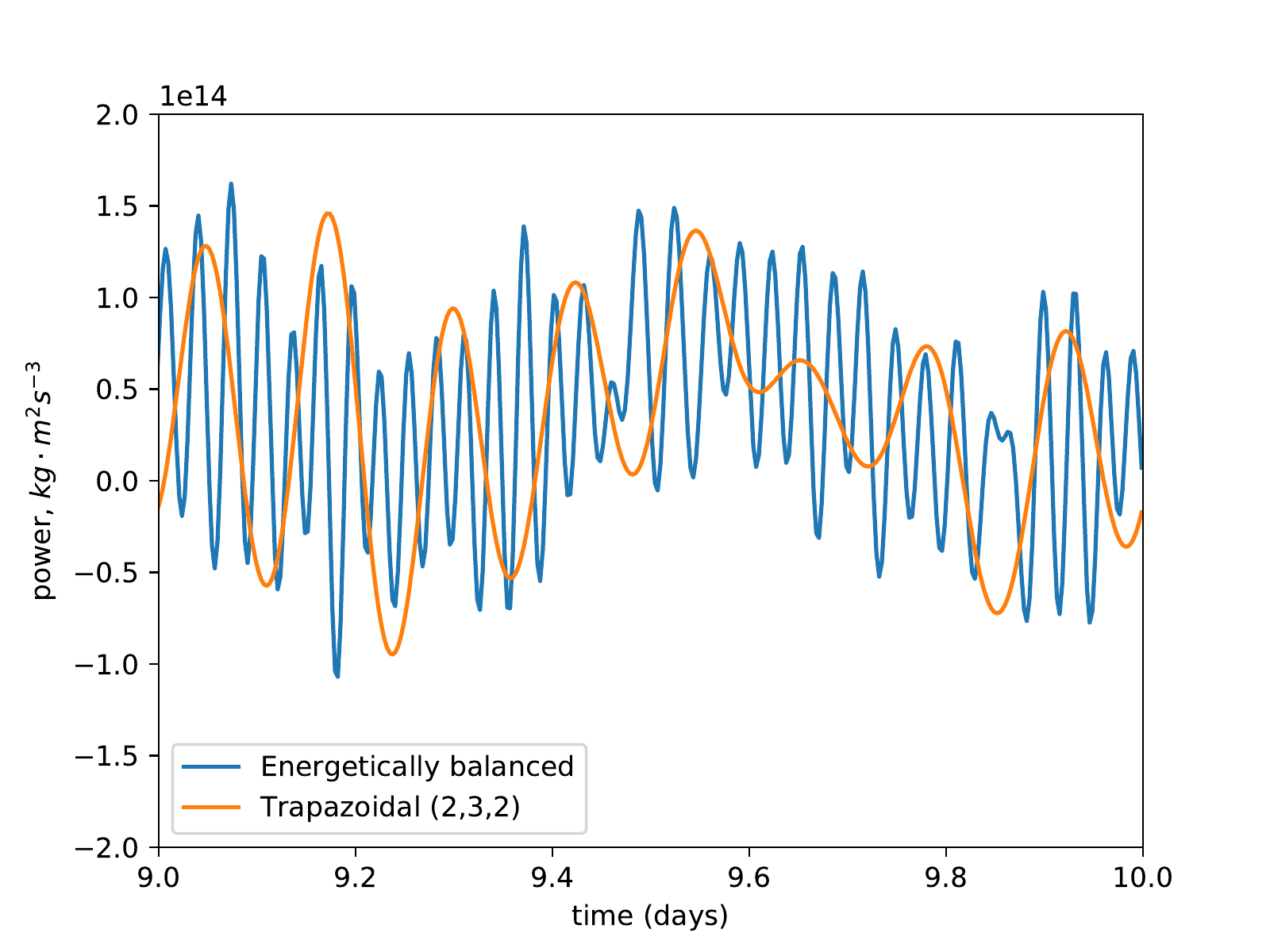}
\includegraphics[width=0.48\textwidth,height=0.36\textwidth]{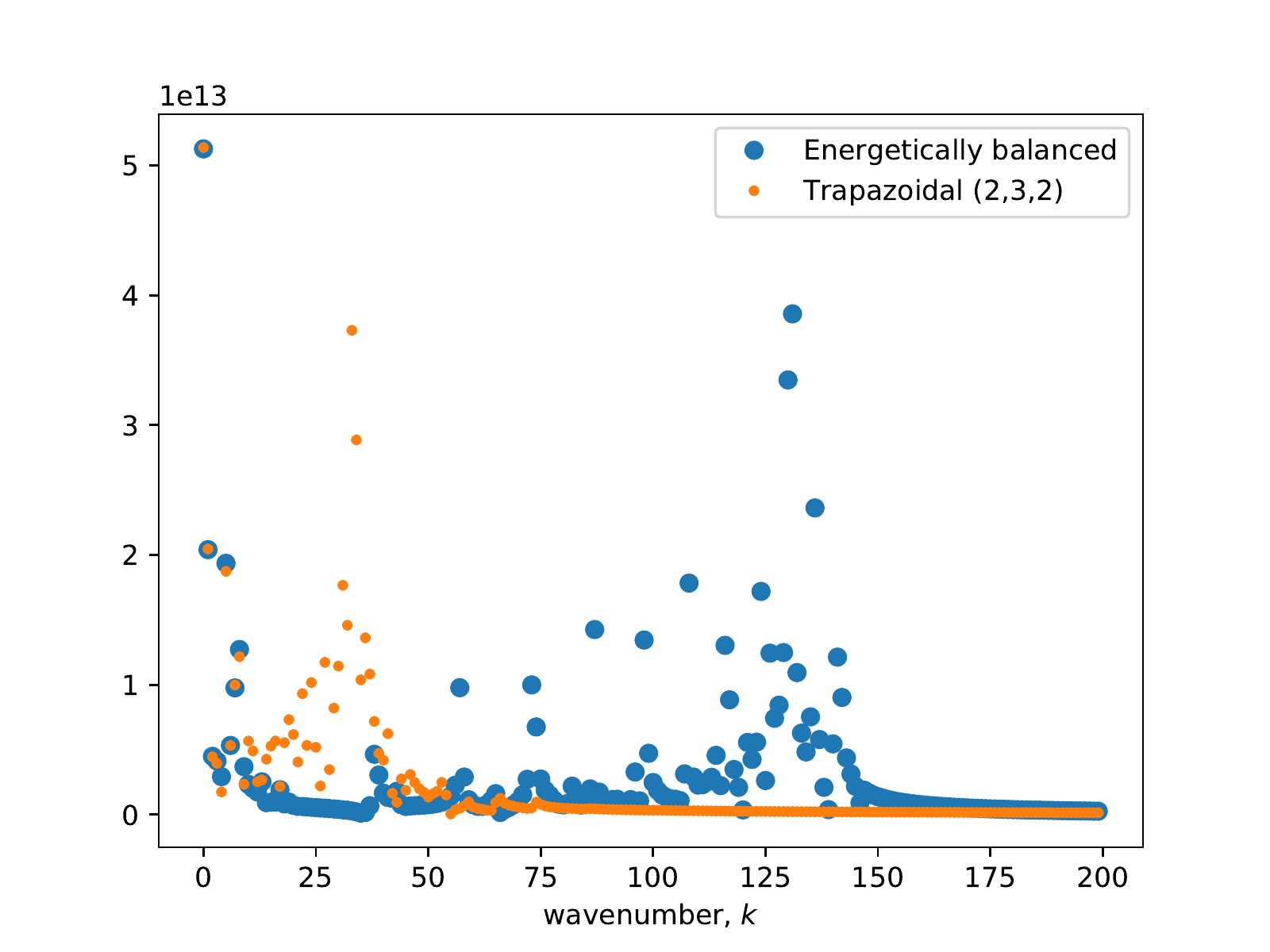}
\caption{\green{Baroclinic instability test case: } Comparison of potential to kinetic \green{energy } 
power exchanges for the energetically balanced and
TRAP(2,3,2) HEVI schemes; left: real space \green{representation } from day 9 to 10, right: amplitude 
of the \green{corresponding } Fourier modes from days 6 to 10.}
\label{fig::p2k}
\end{center}
\end{figure}

The one respect in which the two schemes differ markedly is the potential to kinetic power exchanges,
\blue{the transition of potential to kinetic energy over time. This is computed (assuming homogeneous 
Dirichlet boundary conditions for the vertical velocity) as the discrete equivalent of the term 
$-\partial P/\partial t = g\int z\partial(\rho w)/\partial z\mathrm{d}\Omega$ as
$-\partial P_h/\partial t = g\hat{z}_h\langle\gamma_h,\gamma_h\rangle\boldsymbol{\mathsf{E}}^{3,2}_{\perp}\overline{\hat{W}}_h$
\cite{LP20}. }
For the potential to kinetic power exchange the energetically balanced scheme
exhibits a higher frequency oscillation not present in the TRAP(2,3,2) scheme, as shown in fig. \ref{fig::p2k}.
Fourier analysis of these power exchanges, taken between days 6 and 10, shows that while the amplitudes
of the low-frequency modes ($k<10$) match each other closely, both solutions exhibit a secondary peak,
which for the TRAP(2,3,2) scheme occurs at $k=33$ (2.91 hours), while for the energetically balanced
scheme this peak occurs at $k=131$ (0.73 hours). Curiously, the secondary peak for the energetically
balanced scheme is almost exactly four times the frequency observed for the TRAP(2,3,2) scheme. Since
the secondary peak for the TRAP(2,3,2) scheme is closer in frequency to the low frequency modes, which
exhibit very similar behaviour for both schemes, there is a 
\blue{possibility that this may alias onto the low frequency time scales of the baroclinic motions. }
%danger that this secondary peak
%will alias onto these low frequency modes than there is for the energetically balanced scheme.

\subsection{3D Rising bubble}

In order to verify the scheme for non-hydrostatic dynamics this is applied to a standard test
case for a 3D rising bubble, which is initialised as a small perturbation in potential
temperature against a constant background value, with the Exner pressure specified so as to
satisfy hydrostatic balance \cite{Giraldo13,AG16,Melvin19}. The domain is configured as a
horizontally periodic box of size $1000\mathrm{m}\times 1000\mathrm{m}\times 1500\mathrm{m}$
with a flat bottom and top with homogeneous boundary conditions, and discretised using
$24\times 24$ elements of degree $p=3$ in the horizontal and 150 vertical levels, with a time
step of $\Delta t = 0.01\mathrm{s}$. 
%The model is stabilised by the application of a biharmonic
%viscosity term to the horizontal momentum equation \cite{LP18} with a coefficient of $624.78\mathrm{m}^4\mathrm{s}^{-1}$.
No dissipation of any kind is applied in the vertical, or to the potential temperature advection
equation. \blue{However the energy 
conserving upwinding is applied to the potential temperature diagnostic equation 
\eqref{eq::diag_theta} in order to suppress the development of oscillations associated with
non-hydrostatic motions. This is not applied to the baroclinic instability test case, which 
evolves in a predominantly hydrostatic regime.}

\blue{Since the implicit vertical time stepping scheme conserves energy \cite{Lee20}, and
the energetic exchanges are exactly balanced in the horizontal, no stabilisation is required
of the energetically balanced HEVI scheme, provided that the horizontal CFL condition is satisfied.
This is in contrast to the TRAP(2,3,2) scheme, which requires a horizontal biharminic viscosity in 
order suppress the onset of numerical instability. 
In order to directly compare to the TRAP(2,3,2) scheme, simulations of the energetically balanced
HEVI scheme have been run both with and without horizontal biharmonic viscosity on the momentum
equation \cite{LP18}, with a coefficient of $624.78\mathrm{m}^4\mathrm{s}^{-1}$.}
%Since the implicit vertical time stepping scheme conserves energy \cite{Lee20},
%while the explicit horizontal time stepping scheme presented here does not, stabilisation is
%required only in the horizontal dimensions, even in the presence of strongly oscillatory 
%vertical motions.

\begin{figure}[!hbtp]
\begin{center}
\includegraphics[width=0.48\textwidth,height=0.36\textwidth]{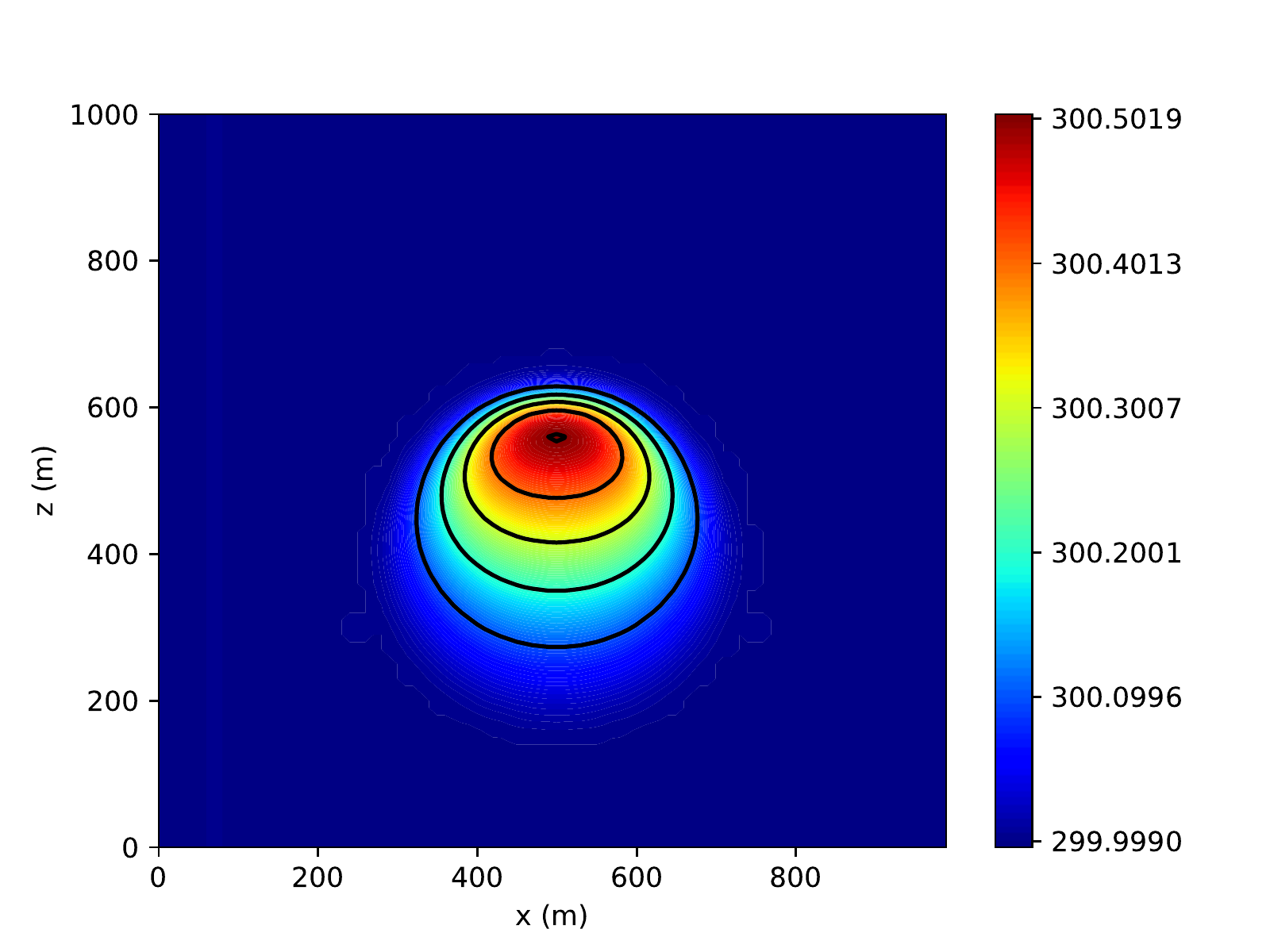}
\includegraphics[width=0.48\textwidth,height=0.36\textwidth]{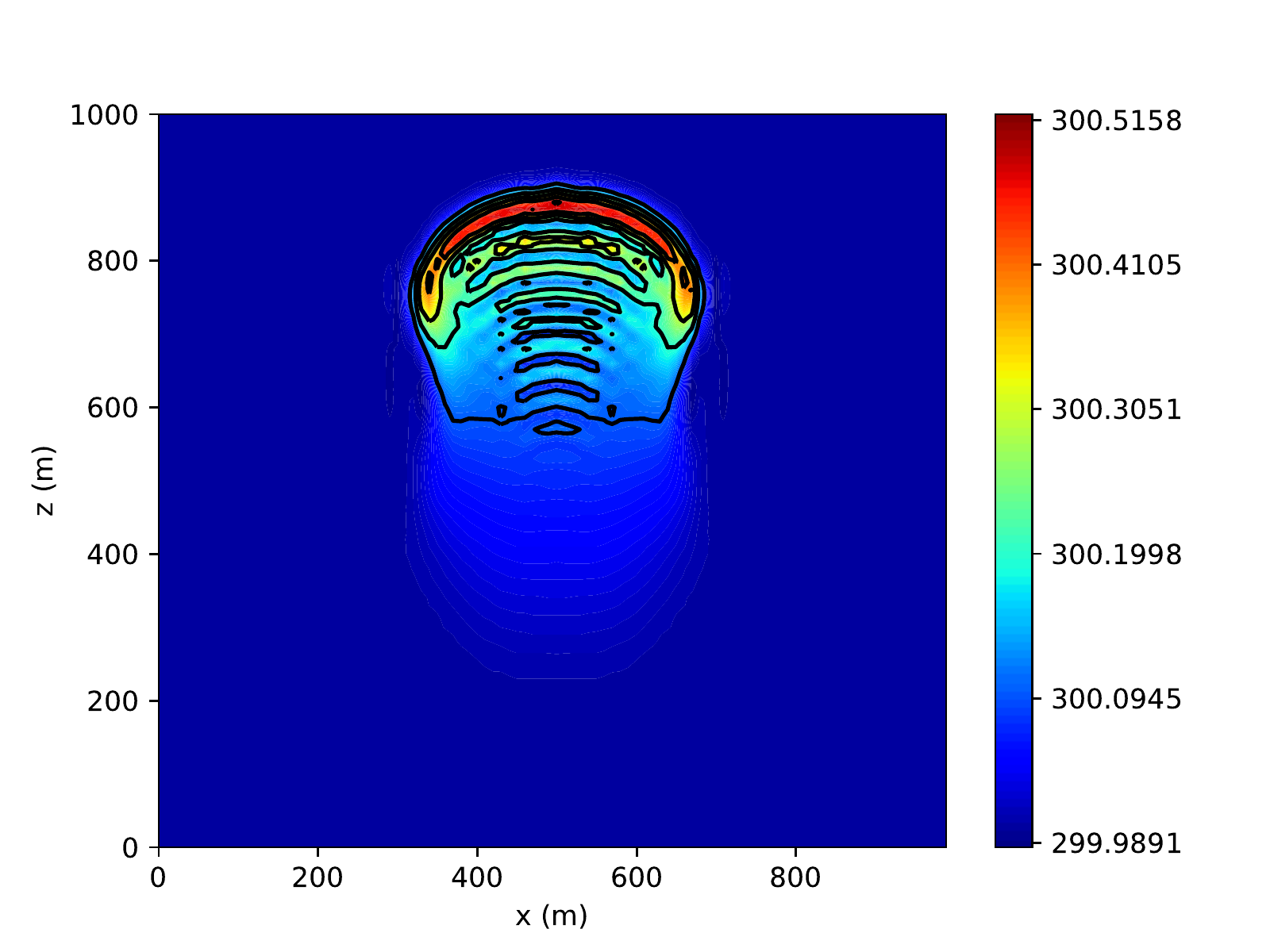}
\caption{\green{3D rising bubble: } Potential temperature cross section at $y=0$ at times 200s (left) and 400s (right).}
\label{fig::bubble}
\end{center}
\end{figure}

\begin{figure}[!hbtp]
\begin{center}
\includegraphics[width=0.48\textwidth,height=0.36\textwidth]{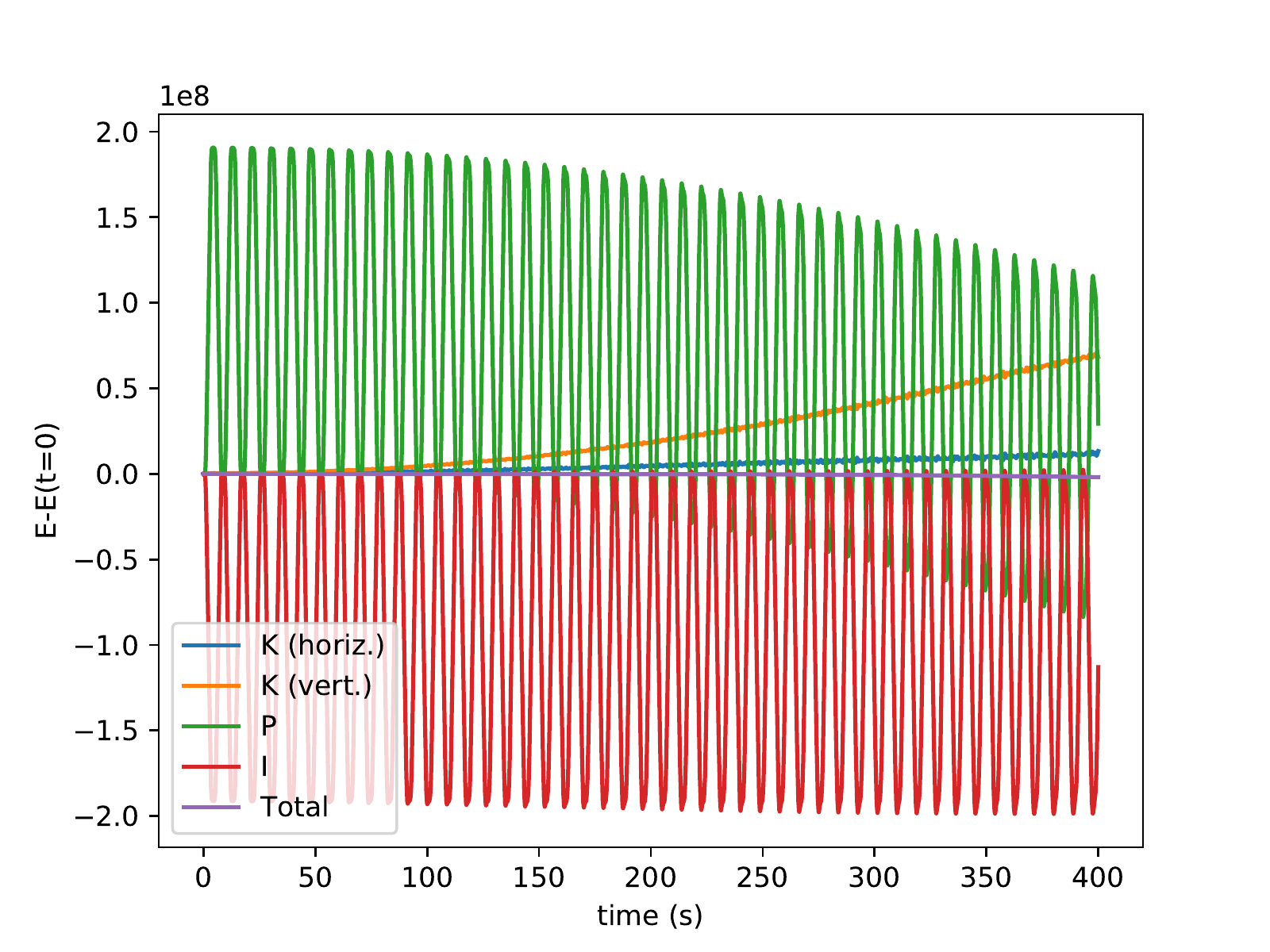}
\includegraphics[width=0.48\textwidth,height=0.36\textwidth]{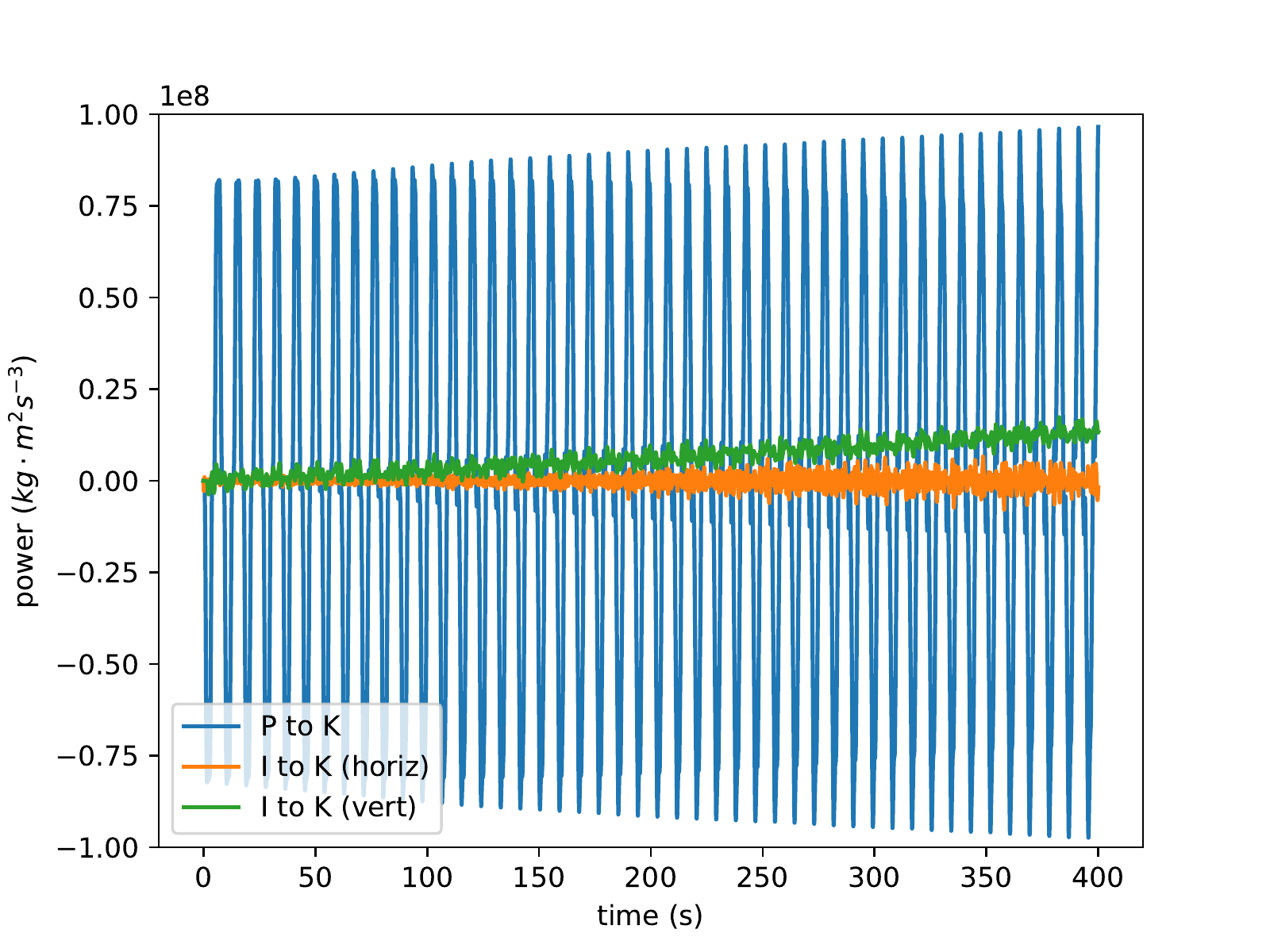}
\caption{\green{3D rising bubble: } Difference in energy from initial value (left), and power exchanges \blue{associated with 
the changes in potential to kinetic and internal to kinetic energy with time } (right). 
\blue{The power exchanges are computed as the discrete analogues of the 
relations $g\int z\partial(\rho w)/\partial z\mathrm{d}\Omega$ and 
$\int\nabla\cdot(\rho\boldsymbol{u}\theta)\Pi\mathrm{d}\Omega$ respectively.}}
\label{fig::bubble_2}
\end{center}
\end{figure}

While the results qualitatively match those previously published in terms of both the
position and deformation of the bubble \cite{AG16,Melvin19}, in the absence of potential temperature
upwinding a secondary oscillation develops
in the wake of the bubble as it steepens, as shown in fig. \ref{fig::bubble}. 
\blue{This oscillation is also observed for the TRAP (2,3,2) scheme \cite{Lee20}. }
The evolution of the energy and power exchanges \blue{(computed as in \cite{LP20}) } are shown in fig 
\ref{fig::bubble_2}
%. These match closely with previous results using a variation on the more 
%common TRAP(2,3,2) HEVI scheme \cite{Lee20}
, and clearly show the signature of the high frequency 
internal gravity wave, as well as the increase in vertical kinetic energy and decrease in 
potential energy which represent the ascension of the bubble.

\subsubsection{Upwinding of potential temperature}

\begin{figure}[!hbtp]
\begin{center}
\includegraphics[width=0.48\textwidth,height=0.36\textwidth]{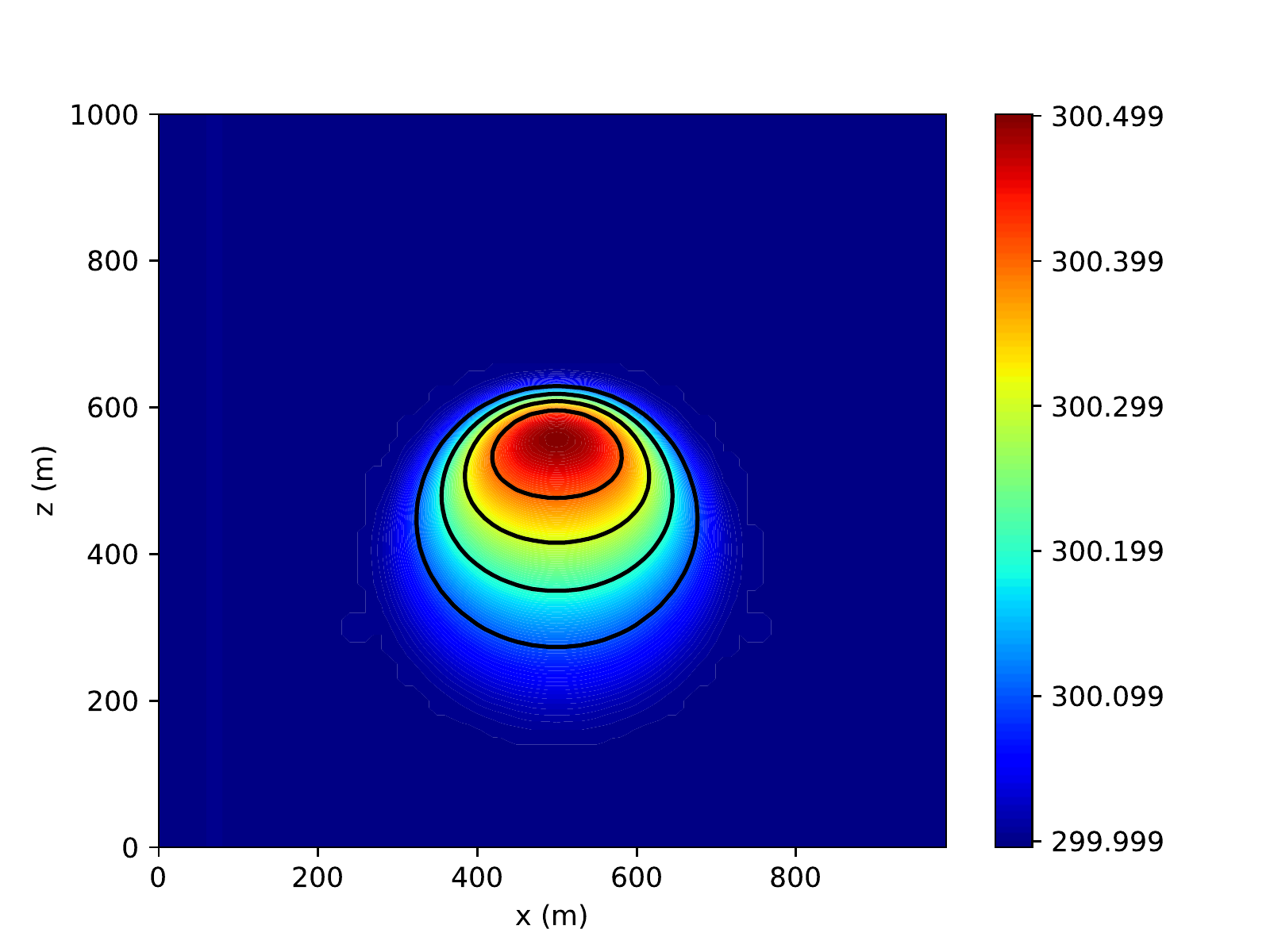}
\includegraphics[width=0.48\textwidth,height=0.36\textwidth]{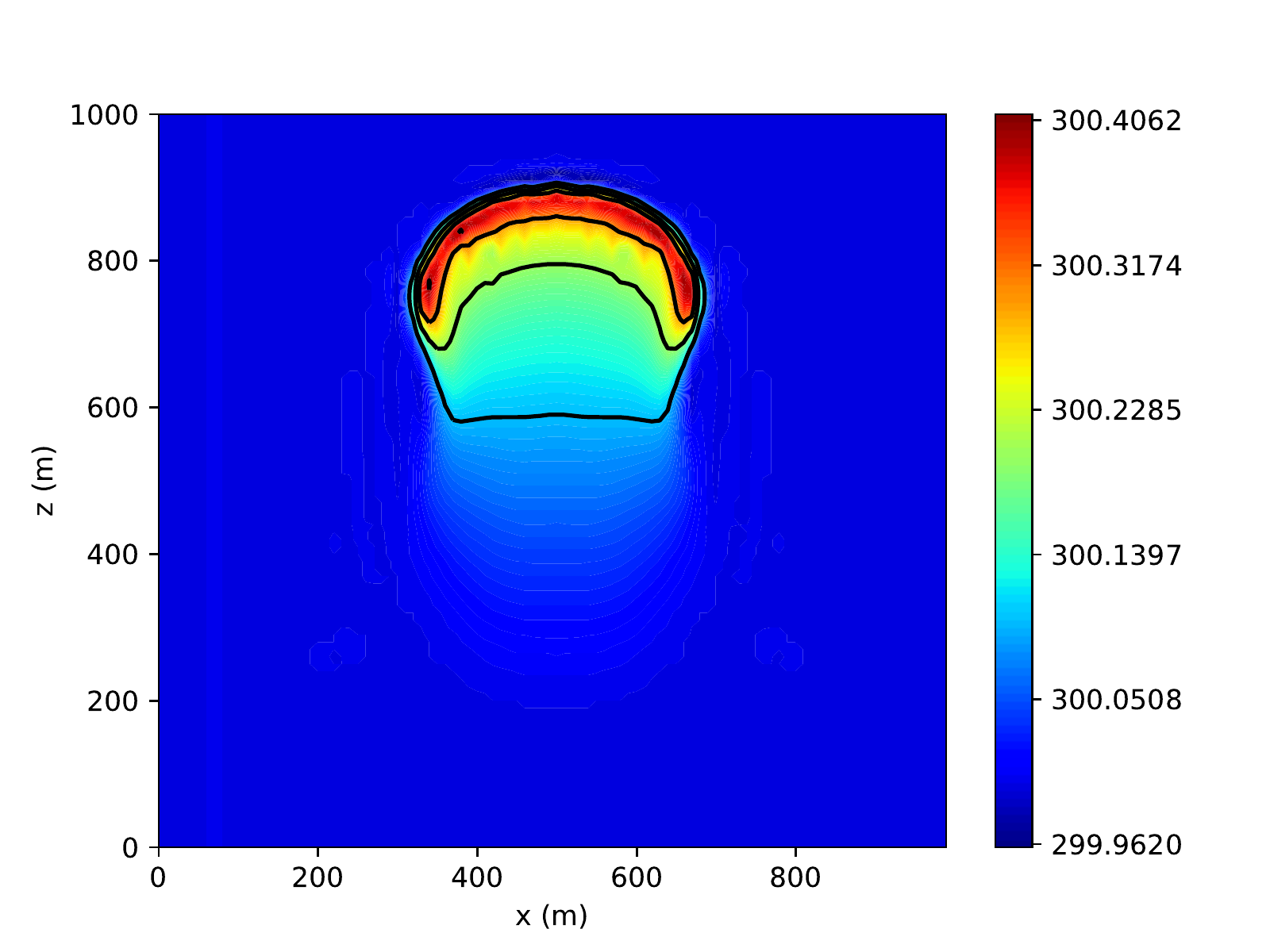}
\caption{\green{3D rising bubble: } Upwinded potential temperature cross section at $y=0$ at times 200s (left) and 400s (right).}
\label{fig::bubble_3}
\end{center}
\end{figure}

\begin{figure}[!hbtp]
\begin{center}
\includegraphics[width=0.48\textwidth,height=0.36\textwidth]{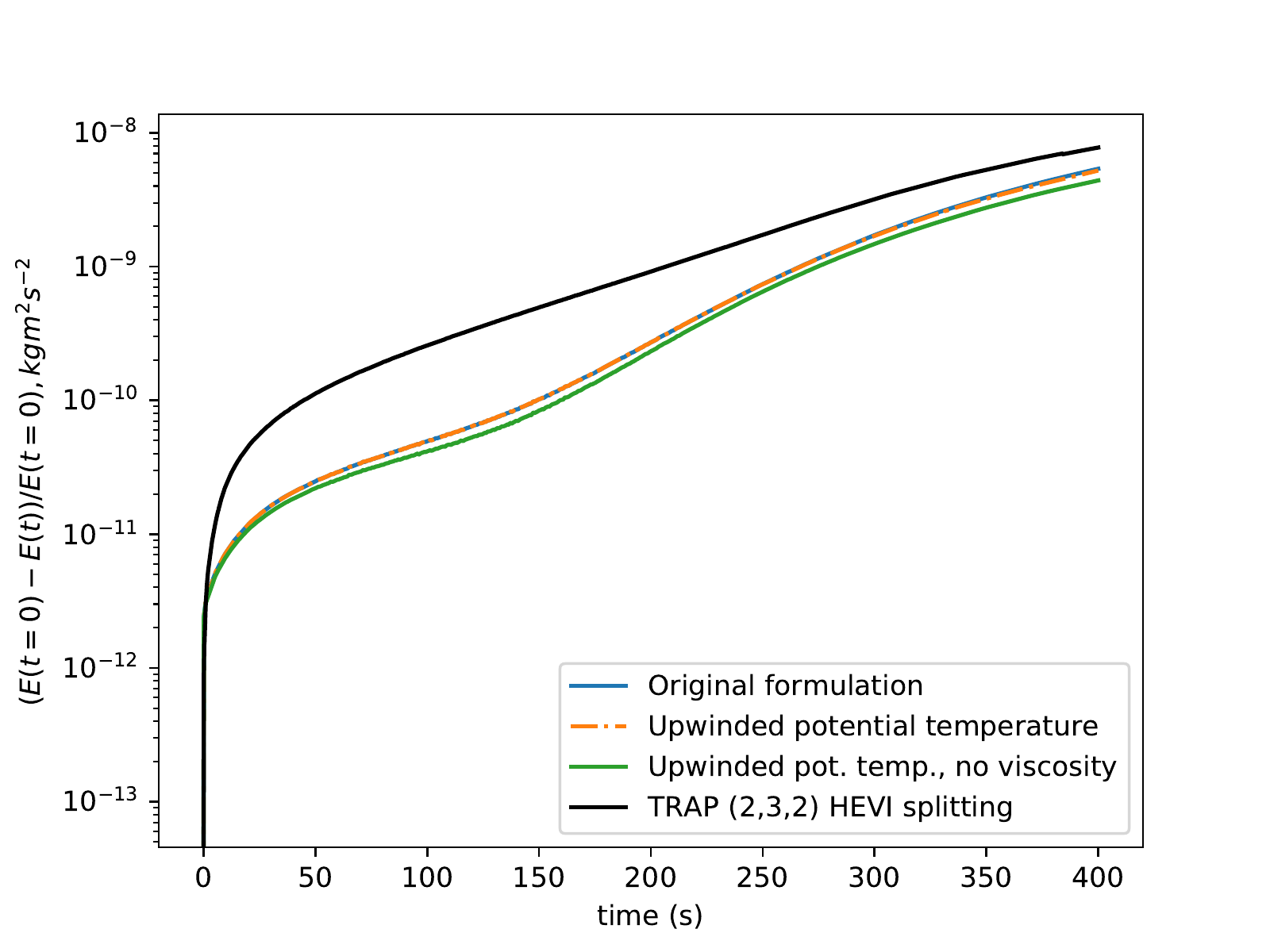}
\includegraphics[width=0.48\textwidth,height=0.36\textwidth]{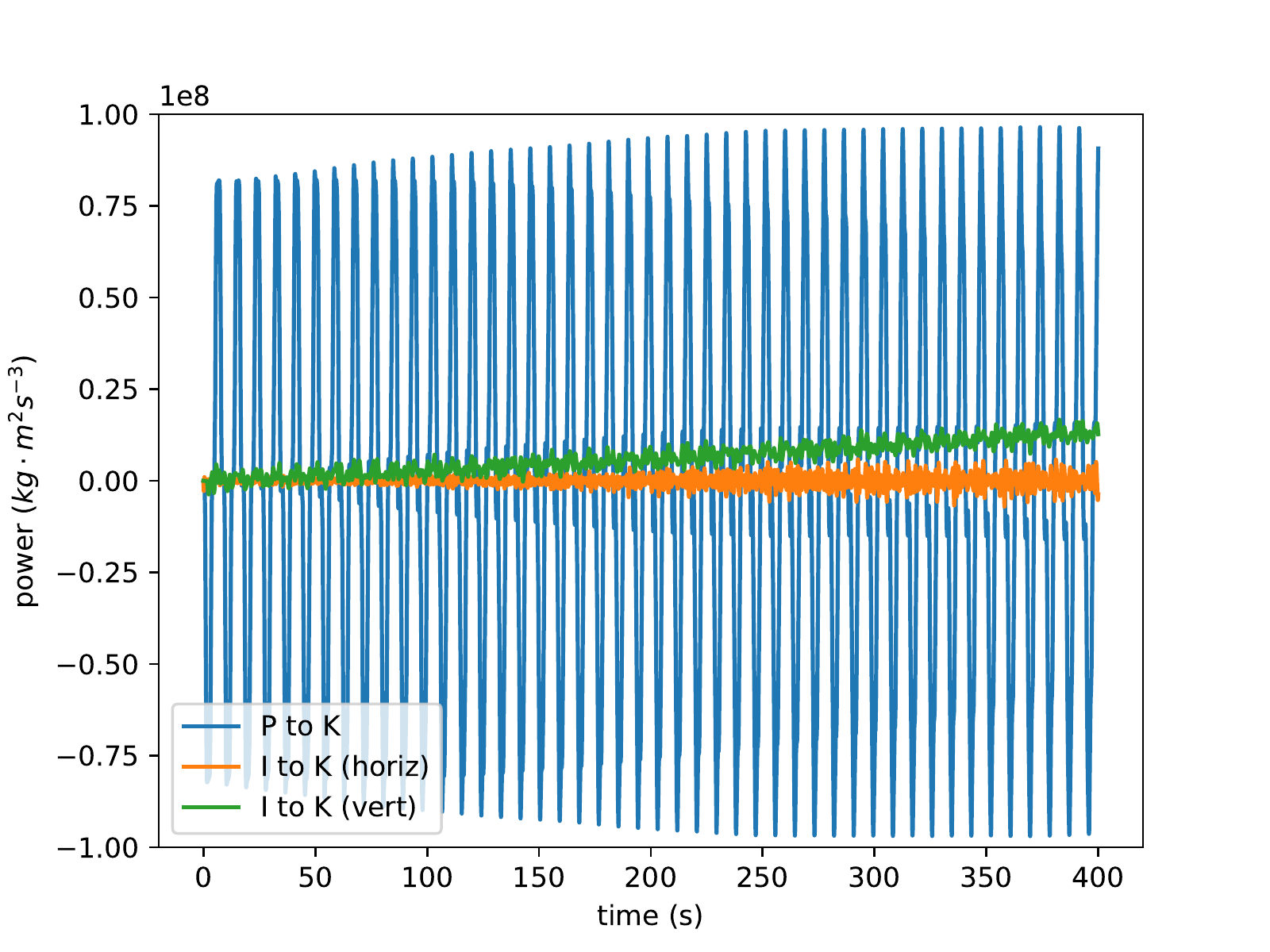}
\caption{\green{3D rising bubble: } Comparison of energy conservation errors for the original and upwinded potential temperature
formulations, as well as the upwinded inviscid solution and the TRAP(2,3,2) scheme (left). \blue{Power exchanges associated with
the changes in potential to kinetic and internal to kinetic energy with time for the TRAP (2,3,2) scheme (right).}}
\label{fig::bubble_4}
\end{center}
\end{figure}

In order to suppress the spurious potential temperature oscillations observed in fig. \ref{fig::bubble}
in an energetically consistent manner, the upwinding formulation presented in section \ref{sec::upwind}
is applied to the rising bubble test case. As observed in fig. \ref{fig::bubble_3}, this upwinding
formulation effectively clears up the spurious overshoots and oscillations. Since the upwinding is performed
within the skew-symmetric operator, as detailed in \eqref{eq::de}, the energy conservation errors are almost
identical to those for the original formulation, as shown in fig. \ref{fig::bubble_4}. 

\blue{Fig. \ref{fig::bubble_4} also shows the energy conservation errors for the TRAP(2,3,2) scheme (with biharmonic viscosity) 
and the energetically balanced scheme in the absence of viscosity. While horizontal biharmonic viscosity is necessary
to stabilise the TRAP(2,3,2) scheme, the energetically balanced scheme may be run entirely without dissipation.
As observed, the difference in the energy conservation error of the energetically balanced scheme with and 
without biharmonic viscosity is small compared to the overall energy conservation error, indicating that most of the
energy conservation error arises from the explicit time stepping of the horizontal velocity, and not the
dissipation term. In all cases the total energy of the system decreases with time.

This figure also shows the power exchanges for the TRAP(2,3,2) scheme. These are almost 
identical to those for the energetically balanced scheme shown in fig. \ref{fig::bubble_2}, such that the 
difference in the time scale of the internal gravity mode observed for the baroclinic test case is not 
observed for the high resolution rising bubble test case, for which the ratio of the horizontal and vertical
resolutions is closer to unity, and the vertical time scales are effectively resolved by the time step.}
%Note that in each case
%a biharmonic viscosity is applied to the horizontal momentum equation, to account for the temporal inconsistency
%in the energetics due to the explicit form of the horizontal time stepping, so exact energy conservation is not
%expected.

\section{Conclusions}\label{sec::conc}

A new horizontally explicit/vertically implicit time splitting scheme for non-hydrostatic atmospheric dynamics
is introduced. The scheme allows for the exact balance of energetic exchanges when coupled with a spatial and
vertical time integration scheme that preserve the skew-symmetric property of the Hamiltonian formulation of
the equations of motion. However since the choice of horizontal velocity for which exchanges balance differs
from the horizontal velocity at the end of the time step, the scheme still permits an energy conservation error.
Linear eigenvalue analysis shows that the second order scheme improves upon the existing second order
trapezoidal HEVI splitting by ensuring the neutral stability of all buoyancy modes, and all acoustic modes below
a certain horizontal CFL number. Comparison to a horizontally third order, vertically second order trapezoidal
splitting for the full 3D compressible Euler equations on the sphere shows that both schemes exhibit a secondary
oscillation in the potential to kinetic power exchanges, associated with an internal gravity mode, however this
oscillation has a temporal frequency approximately four times faster for the new energetically balances scheme,
\blue{such that the time scales of the oscillation for the trapezoidal scheme are closer to those of the baroclinic motions.}
%meaning this mode is less likely to alias onto the baroclinic time scales than in the traditional scheme.

An energetically consistent formulation for potential temperature upwinding is also presented, by which the
test functions are evaluated at vertically downstream locations. This ensures that spurious oscillations in
potential temperature can be removed for vertical motions at non-hydrostatic regimes. This upwinding is
applied within the skew-symmetric formulation, so as to suppress oscillations associated with vertical motions
without altering the energetic properties of the model, for which the exchanges remain balanced.

\blue{Notably, for the rising bubble test case, which is run on an affine geometry, the energetically balanced 
HEVI scheme may be run stably without dissipation of any kind. For the baroclinic instability test case on the
sphere however, horizontal biharmonic viscosity is still required in order to stabilise the simulation. 
Since the model uses inexact integration, and the sphere is represented as a non-affine geometry, it is possible
that aliasing errors may provide an additional source of instability. As part of future work exact integration
will be employed in order to determine if the scheme can be run entirely without dissipation for the
baroclinic instability test case on the sphere also.}

The code used in this article is publicly available at: \verb|https://github.com/davelee2804/MiMSEM|.

\section*{Appendix A: Temporal discreisation of the energetically balanced HEVI scheme}

\blue{For the benefit of readers who are not familiar with compatible finite element methods, in this 
appendix we present a formulation of the energetically balanced HEVI scheme that is discretised in the 
temporal dimension only. Note that in order to preserve energetic balance a spatial discretisation that
respects the skew-symmetric structure of the equations of motion must also be employed. The temporal 
scheme is expressed as follows:}

\emph{Step 1: Explicit horizontal momentum solve}

\begin{equation}
\boldsymbol{v}' = \boldsymbol{v}^n - \Delta t\ {q}_{\perp}^{n}\times\boldsymbol{V}^n
+\Delta t\ \boldsymbol{q}_{\parallel}^{n}W^n
-\Delta t\ \nabla_h\Phi^n - \Delta t\ \theta^{n}\nabla_h\Pi^n
\end{equation}

\emph{Step 2: Implicit vertical solve (including horizontal divergence)}

\begin{subequations}
\begin{align}
w^{n+1} &= w^n - \Delta t\ \partial_z\overline{\Phi} - \Delta t\ \theta^{n+1/2}\partial_z\overline{\Pi} -
\Delta t\ \boldsymbol{q}_{\parallel}^{n+1/2}\cdot\overline{\boldsymbol{V}} \\
\rho^{n+1} &= \rho^n - \Delta t\ \partial_z\overline{W} - \Delta t\ \nabla_h\cdot\overline{\boldsymbol{V}}\\
\Theta^{n+1} &= \rho^n - \Delta t\ \partial_z(\theta^{n+1/2}\overline{W}) - \Delta t\ \nabla_h\cdot(\theta^{n+1/2}\overline{\boldsymbol{V}})
\end{align}
\end{subequations}
with
\begin{subequations}
\begin{align}
\overline{V} &= \frac{1}{6}\Big(2\rho^{n}\boldsymbol{v}^{n} + \rho^{n+1}\boldsymbol{v}^{n} + 
\rho^{n}\boldsymbol{v}' + 2\rho^{n+1}\boldsymbol{v}'\Big) \\
\overline{W} &= \frac{1}{6}\Big(2\rho^{n}w^{n} + \rho^{n+1}w^{n} + \rho^{n}w^{n+1} + 2\rho^{n+1}w^{n+1}\Big) \\
\overline{\Phi} &= \frac{1}{6}\Big(\boldsymbol{v}^{n}\boldsymbol{v}^{n} + \boldsymbol{v}^{n}\boldsymbol{v}' + \boldsymbol{v}'\boldsymbol{v}' + 
w^{n}w^{n} + w^{n}w^{n+1} + w^{n+1}w^{n+1}\Big) + gz \\
\overline{\Pi} &= \frac{1}{2}\Big(\Pi^{n} + \Pi^{n+1}\Big)
\end{align}
\end{subequations}

\emph{Step 3: Explicit horizontal momentum solve}

\begin{equation}
\boldsymbol{v}^{n+1} = \boldsymbol{v}^n - \Delta t\ {q}_{\perp}^{n+1/2}\times\overline{\boldsymbol{V}}
+\Delta t\ \boldsymbol{q}_{\parallel}^{n}\overline{W}
-\Delta t\ \nabla_h\overline{\Phi} - \Delta t\ \theta^{n+1/2}\nabla_h\overline{\Pi}
\end{equation}
where $\nabla_h$ is the horizontal gradient operator.

Note that if energetic balance is not of concern, a standard time centered construction of the mass flux 
terms, $\boldsymbol{V}$, $W$ and Bernoulli function, $\Phi$ at $n+1/2$ should also preserve second order 
accuracy. Also note that a simple first order Euler integration is used here in step 1, which is sufficient 
to construct a second order mass flux and Bernoulli function.

\section*{Appendix B: Entropy conservation for the variational form of the compressible Euler equations}

The time derivative of the internal energy is derived
by multiplying the density-weighted potential temperature advection equation,
\begin{equation}
\frac{\partial\rho\theta}{\partial t} + \nabla\cdot(\rho \boldsymbol{u}\theta) = 0,
\end{equation}
by the Exner pressure,
\begin{equation}
\frac{\delta I}{\delta\Theta} = \Pi = c_p\Bigg(\frac{R\rho\theta}{p_0}\Bigg)^{R/c_v},
\end{equation}
such that the internal energy, ${I}$ evolves as
\begin{equation}
\frac{\partial{I}}{\partial t} =
\Bigg\langle\frac{\delta I}{\delta\Theta},\frac{\partial\Theta}{\partial t}\Bigg\rangle =
\Bigg\langle\Pi,\frac{\partial\rho\theta}{\partial t}\Bigg\rangle =
\frac{c_v}{c_p}\frac{\partial}{\partial t}\int\rho\theta\Pi\mathrm{d}\Omega =
-\langle\Pi,\nabla\cdot(\boldsymbol{u}\rho\theta)\rangle,
\end{equation}
since $R = c_p - c_v$ and $\Theta = \rho\theta$. This expression may then be expanded as:
\begin{subequations}
\begin{align}
\frac{c_v}{c_p}\frac{\partial}{\partial t}\int\rho\theta\Pi\mathrm{d}\Omega
&= -\langle\Pi,\theta\nabla\cdot(\rho \boldsymbol{u}) + \rho \boldsymbol{u}\cdot\nabla\theta\rangle \\
&= -\Bigg\langle\rho\theta\Pi,\frac{1}{\rho}\nabla\cdot(\rho \boldsymbol{u}) + \frac{\boldsymbol{u}}{\theta}\cdot\nabla\theta\Bigg\rangle \\
&= \Bigg\langle\rho\theta\Pi,\frac{\partial\log(\rho)}{\partial t} - \boldsymbol{u}\cdot\nabla\log(\theta)\Bigg\rangle
\end{align}
\end{subequations}
The left hand side of this expression may be re-arranged as:
\begin{equation}
\frac{c_v}{c_p}\int\frac{\partial\rho\theta\Pi}{\partial t}\mathrm{d}\Omega =
\frac{c_v}{c_p}\Bigg\langle\frac{\rho\theta\Pi}{\rho\theta\Pi},\frac{\partial\rho\theta\Pi}{\partial t}\Bigg\rangle =
\frac{c_v}{c_p}\Bigg\langle\rho\theta\Pi,\frac{\partial\log(\rho\theta\Pi)}{\partial t}\Bigg\rangle =
\Bigg\langle\rho\theta\Pi,\frac{\partial\log(\rho)}{\partial t} - \boldsymbol{u}\cdot\nabla\log(\theta)\Bigg\rangle
\end{equation}
For $\rho\theta\Pi\in L^2(\Omega)$, cancellation of this factor from both sides yields the variational
expression
\begin{equation}
\frac{c_v}{c_p}\Bigg\langle\gamma,\frac{\partial\log(\rho\theta\Pi)}{\partial t}\Bigg\rangle =
\Bigg\langle\gamma,\frac{\partial\log(\rho)}{\partial t} - \boldsymbol{u}\cdot\nabla\log(\theta)\Bigg\rangle\qquad\forall\gamma\in L^2(\Omega)
\end{equation}
Recalling the equation of state gives
\begin{equation}
\rho\theta\Pi = c_p\Bigg(\frac{R}{p_0}\Bigg)^{R/c_v}(\rho\theta)^{c_p/c_v},
\end{equation}
\begin{equation}
\Bigg\langle\gamma,\frac{\partial\log(\rho\theta)}{\partial t}\Bigg\rangle =
\Bigg\langle\gamma,\frac{\partial\log(\rho)}{\partial t} - \boldsymbol{u}\cdot\nabla\log(\theta)\Bigg\rangle\qquad\forall\gamma\in L^2(\Omega)
\end{equation}
which is further simplified as:
\begin{equation}
\Bigg\langle\gamma,\frac{\partial\log(\theta)}{\partial t} + \boldsymbol{u}\cdot\nabla\log(\theta)\Bigg\rangle = 0\qquad\forall\gamma\in L^2(\Omega)
\end{equation}
Recalling the original definition of the entropy then gives
\begin{equation}
\Bigg\langle\gamma,\frac{\partial s}{\partial t} + \boldsymbol{u}\cdot\nabla s\Bigg\rangle = 0\qquad\forall\gamma\in L^2(\Omega)
\end{equation}
such that entropy is materially conserved. Setting the test function as the density $\gamma = \rho$ in the above expression,
and multiplying the continuity equation by the entropy, we have
\begin{subequations}
\begin{align}
\Bigg\langle\rho,\frac{\partial s}{\partial t} + \boldsymbol{u}\cdot\nabla s\Bigg\rangle &= 0\\
\Bigg\langle s,\frac{\partial\rho}{\partial t} + \nabla\cdot\boldsymbol{u}\rho\Bigg\rangle &= 0.
\end{align}
\end{subequations}
%Adding the above expressions gives
Adding this to the continuity equation, and recalling that $\boldsymbol{U} := \rho\boldsymbol{u}$ gives
\begin{equation}
\int\frac{\partial\rho s}{\partial t} + \nabla\cdot(\boldsymbol{u}\rho s)\mathrm{d}\Omega = 0,
\end{equation}
which yields an additional conservation law for the entropy function-entropy flux pair $\rho s, \boldsymbol{u}\rho s$ \cite{Tadmor13}.

\section{Acknowledgments}

David Lee would like to thank Dr Justin Freeman for his continued encouragement, and 
Drs. Marcus Thatcher and John McGregor for their support and access to computing resources.
This project was supported by resources and expertise provided by CSIRO IMT Scientific Computing.
We are also grateful to the two anonymous reviewers, whose comments and insights helped to 
improve the quality of this article.

\end{document}